\documentclass[10pt,a4paper,reqno]{article}

\usepackage{amsfonts,amsthm,amssymb,amsmath,graphicx,enumerate}
\usepackage{mathrsfs,setspace,multirow, multicol,latexsym}
\usepackage{mathtools}
\usepackage{tikz-cd}
\usepackage{courier}
\usepackage{amsfonts}%math fonts
\usepackage{multicol}%for multiple columns
\usepackage[mathscr]{euscript}%for script letters
\usepackage{dirtytalk}
\usepackage{enumitem}
\usepackage{dirtytalk}
\usepackage{hyperref}
\usepackage{pdfpages}

\usepackage{graphicx}
\usepackage[list=true,listformat=simple,labelformat=simple]{subcaption}

\newtheorem{theorem}{{Theorem}}[section]

\newtheorem{lemma}[theorem]{{Lemma}}
\theoremstyle{definition}
\newtheorem{definition}[theorem]{{Definition}}%[section]

\newtheorem{example}[theorem]{{Example}}
\newtheorem{notation}[theorem]{{Notation}}

\theoremstyle{remark}
\newtheorem{remark}[theorem]{{Remark}}

\newcommand{\R}{\mathbb R}

\newcommand{\N}{\mathrm N}

\newcommand{\M}{\mathcal M}
\newcommand{\A}{\mathcal A}
\newcommand{\I}{\mathcal I}

\newcommand{\LP}{\left(}
\newcommand{\RP}{\right)}
\newcommand{\LS}{\left[}
\newcommand{\RS}{\right]}
\newcommand{\LC}{\left\{}
\newcommand{\RC}{\right\}}
\newcommand{\MV}{\,\middle \vert\,}

\title{Polynomially knotted $2$ Spheres}
\author{Tumpa Mahato\thanks{Department of Mathematics, Indian Institute of Science Education and Research, Pune, India. Email: tumpa.mahato@students.iiserpune.ac.in} and Rama Mishra \thanks{Department of Mathematics, Indian Institute of Science Education and Research, Pune, India. Email: r.mishra@iiserpune.ac.in}}

\date{\today}

\begin{document}

\maketitle

\begin{abstract}
We review the polynomial parameterization of classical knots and prove the analogous results for long $2$ knots. We also construct polynomial parameterization for certain classes of knotted spheres (such as spun and twist spun of the classical knots).

\end{abstract}

\section{Introduction} 
Higher dimensional knot theory has attracted a lot of interest among knot theorists recently. By an $n$-knot, one means a smooth, proper, locally flat embedding of $S^n$ in $S^{n+2}$. Moving from classical knots, the simplest situation is to understand knotting of $S^2$ in $S^4$ or in $\R^4$. In fact, studies are being done on any surface getting knotted inside four-dimensional space, and $2$-knots are just a special case. The notion of ambient isotopy can be defined
as in the case of classical knots. The central problem remains the same: classify all $2$-knots up to ambient isotopy.

Similar to the classical knots, a diagrammatic theory (\cite{CarSai}) and a braid theory (\cite{Kam}) have been developed for $2$-knots. Both theories have provided many interesting invariants for $2$-knots. In the case of classical knots, we have encountered many numerical invariants whose computations become easy using a suitable parameterization.  In knot theory, parameterization of knots is a way to present a knot, defining its coordinates using elementary functions, such as trigonometric or polynomial functions. They are referred to as trigonometric knots or polynomial knots.
Fourier knots \cite{Kau98} and Lissajous knots \cite{Bog} are examples of trigonometric knots. Trigonometric parameterization provides a compact classical knot embedding, whereas polynomial parameterization presents long knots, i.e., embeddings of $\R$ in $\R^{3}$. Polynomial knots are important because they bring an algebraic flavour into knot theory. One point compactification of a polynomial knot is a classical knot, and the projective closure of a polynomial knot is a knot in $\R P^3.$  Polynomial knots are extensively studied by Shastri \cite{Sha}, Mishra-Prabhakar \cite{PraMis}, Durfee-Oshea \cite{durfoshea}.
Explicit polynomial representation for all knots up to $8$ crossings can be found in \cite{PraMis}.

In this paper, we discuss parameterization of $2$-knots using polynomial functions. Since a $2$-knot is a proper, smooth, locally flat embedding of $S^{2}$ in $S^4$, removing a point from $S^2$ and its corresponding image under the embedding in $S^{4}$ we get an embedding of $\R^{2}$ to $\R^{4}$. We will call this embedding a \textit{long $2$-knot}. We prove that {\it every long $2$-knot is ambient isotopic to a polynomial embedding of $\R^2$ in $\R^{4}$}. We also show that any two polynomial embeddings of $\R^2$ in $\R^4$, which are ambient isotopic, are, in fact, polynomially isotopic. The idea of our proof is similar to the classical case \cite{Shu}.

We also address the problem of explicitly constructing a polynomial parameterization for a given $2$-knot. We begin by taking the simplest $2$-knots, which are constructed from classical knots, such as \textit{spun knots}, introduced by E. Artin \cite{Art}. Then we move to {\it d-twist spun 2-knots } which is a generalization of Artin's construction, introduced by E. C. Zeeman \cite{Zee65}.  Obtaining a polynomial parameterization for  
these classes of knots is relatively simple because we can use the existing polynomial representation of classical knots \cite{PraMis}
which can be brought into the form suitable for spun or twist spun construction. In this paper, we provide the steps to parameterize spun knots and twist spun knots using polynomial functions and present their $3$ dimensional projections using Mathematica.

This paper is organized as follows: Section \ref{sec:pre} is divided into three subsections. Section \ref{sec:surface} includes the basic definitions required to study $2$-knots. Section \ref{sec:sphere} describes the two constructions, spinning and twist spinning, for obtaining 2-knots from classical knots. In Section \ref{sec:polyknots}, we provide a brief exposition of the polynomial parameterization of classical knots. In Section \ref{sec:long}, we prove the existence of polynomial parameterization for all long $2$-knots and show that it is unique up to a polynomial isotopy.  In Sections \ref{sec:spunPoly}  and \ref{th:twistPoly}, we provide  steps to find a polynomial representation for spun knots and twist spun knots, respectively. We write down explicit parameterizations for the spun trefoil and the spun figure eight knot and also the $d$ twist spun trefoil for $d=1,2,5,10$ and $20$. Finally, in Section \ref{mathematica}, the {\it Mathmatica} notebooks for all the examples are provided.
\section{Prerequisites}\label{sec:pre}

\subsection{Surface knots}\label{sec:surface}

In order to understand surface knots in detail, we present some definitions here.

\begin{figure}[h]
	\centering
	\includegraphics[width=0.5\textwidth]{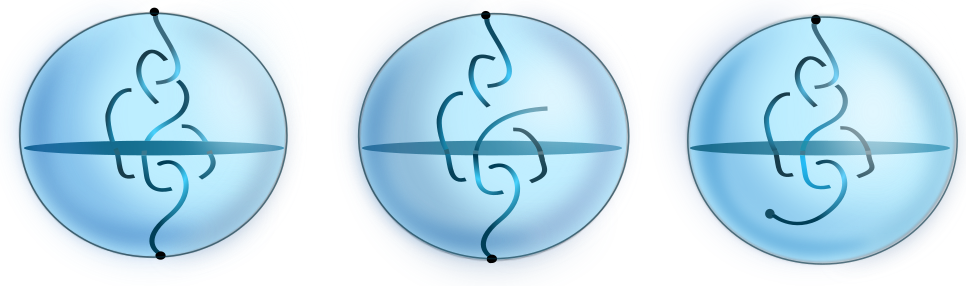}
	\caption{Examples of proper and non-proper embeddings}
	\label{fig: proper}
\end{figure}

\begin{definition}
	
	Let $N$ be an $n$-dimensional manifold embedded in an $m$-dimensional manifold $M$. We call an embedding $f:N \rightarrow M$ to be {\it proper} if, every compact subset $K \subset f(N)$ has a compact preimage $f^{-1}(K)=\{x \in N \vert f(x) \in K\}$.
	
	In other words, the manifold $f(N)$ is proper if $\partial f(N)= f(N) \cap \partial M$.
\end{definition}
\begin{example}
	The left most figure of Figure \ref{fig: proper} shows an example of a properly embedded 1-manifold in a 3-ball. The figure in the middle is not proper because one of the interior points of the 1-manifold is on the boundary of the 3-ball. The right most figure is not proper because one of the boundary points of the $1$-manifold is an interior point of the 3-ball. 
\end{example}
\begin{definition}
	Let $f:N\to M$ be a proper embedding of smooth manifolds. This embedding is said to be 
  {\it  locally flat}  at a point $x \in M$ if there exists a regular neighborhood $U$ of $x$ in $M$ such that the pair $(U,U \cap f(N))$ is homeomorphic to the standard ball pair $(D^{m},D^{n})$, where $D^{k}$ is the unit ball in $\R^{k}$, centred at origin.
	We say that $f$ is {\it locally flat} if it is locally flat at every point of $M$.
\end{definition}

Non-locally flat embeddings exist.  An example can be seen in \cite{Kam}.
\begin{definition}
	A {\it surface link} is a proper, locally flat embedding of a closed surface $F$ in $\R^{4}$ or $S^{4}$. 
\end{definition}
\begin{itemize} [topsep=0pt]
	\item If $F$ is a connected surface then it is called a {\it surface knot}.
	\item When $F=S^{2}$, it is called a {\it 2-knot} or {\it $2$-dimensional knot}.
\end{itemize}
\begin{definition}
	
	A proper, locally flat  embedding of $\R^{2}$  in $\R^{4}$ which is asymptotic outside a compact region is referred to as a  {\it long $2$-knot}.
	
\end{definition}

Note that every  $2$-knot will give rise to long $2$-knot by removing one point from $S^{2}$ and the image of that point from $S^{4}$ and considering the resulting restriction of the embedding.
\begin{definition}
	Two surface links $F$ and $F'$ are {\it equivalent} if and only if they are ambient isotopic, i.e. there exists an ambient isotopy $\{h_{t}\}$ of $\R^{4}$ such that $h_{1}(F)=F'$.
	
	When $F$ and $F'$ are oriented, then $h_{1}\vert_{F}:F\rightarrow F'$ is assumed to be an orientation-preserving homeomorphism.
\end{definition}

\begin{figure}[h]
	\centering
	\begin{subfigure}[c]{0.3\linewidth}
		\centering
		\includegraphics[width=0.6\textwidth]{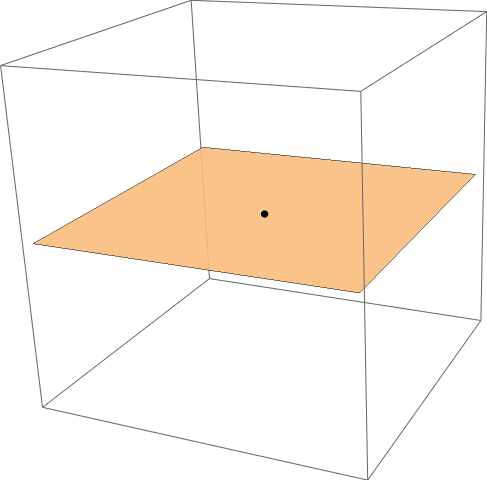}
		\caption{Regular point}
		\label{fig:regularpoint}
	\end{subfigure}
	\quad
	\begin{subfigure}[c]{0.3\linewidth}
		\centering
		\includegraphics[width=0.65\textwidth]{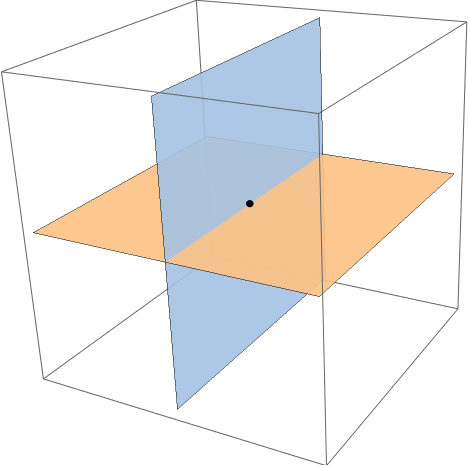}
		\caption{Double point}
		\label{fig:doublepoint}
	\end{subfigure}
	\quad
	\begin{subfigure}[c]{0.3\linewidth}
		\centering
		\includegraphics[width=0.65\textwidth]{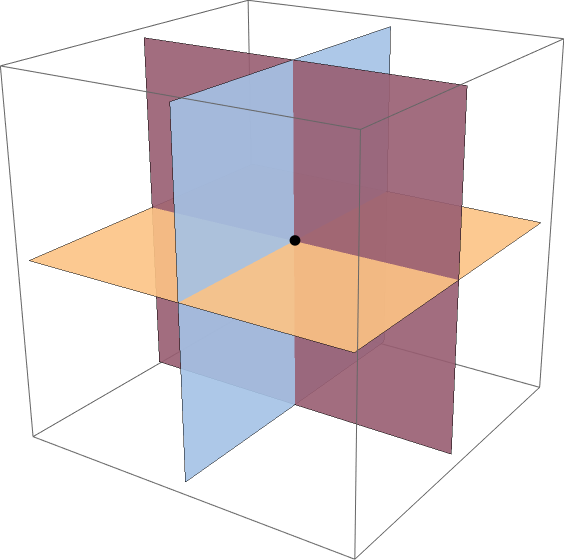}
		\caption{Triple point}
		\label{fig:triplepoint}
	\end{subfigure}
	\newline 
	\begin{subfigure}[c]{0.3\linewidth}
		\centering
		\includegraphics[width=0.7\textwidth]{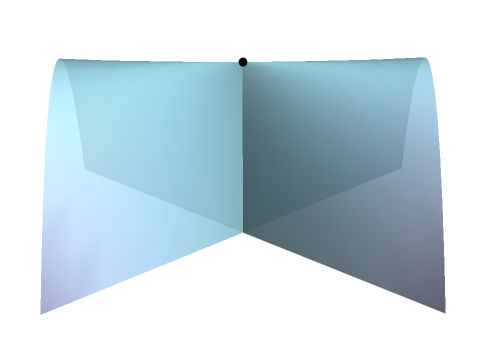}
		\caption{Branch point}
		\label{fig:branchpoint}
	\end{subfigure} 
	\caption{Neighborhood of singular points}
	\label{fig:singular points}
\end{figure}

\begin{figure}[h]
  	\centering
  	\begin{subfigure}[b]{0.3\linewidth}
  		\centering
  		\includegraphics[width=0.68\textwidth]{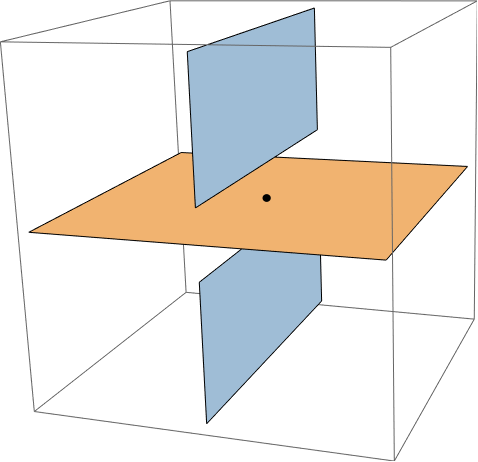}
  		\caption{Near a double point}
  	\end{subfigure}
  	\begin{subfigure}[b]{0.3\linewidth}
  		\centering
  		\includegraphics[width=0.68\textwidth]{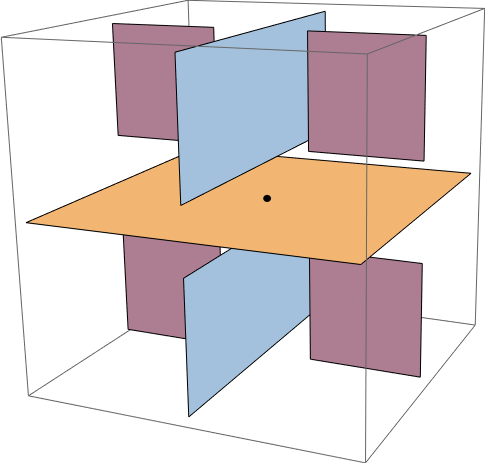}
  		\caption{Near a triple point}
  	\end{subfigure} 
  	\begin{subfigure}[b]{0.3\linewidth}
  		\centering
  		\includegraphics[width=0.68\textwidth]{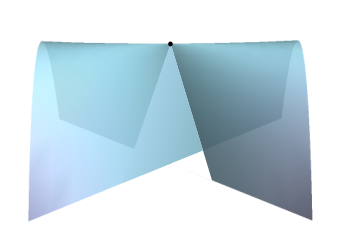}
  		\caption{Near a branch point}
  	\end{subfigure} 
  	\caption{Broken surface diagram near a singular points}
  	\label{fig:broken}
  \end{figure}

\begin{figure}[h]
  	\centering
  	\begin{subfigure}[b]{0.45\linewidth}
  		\centering
  		\includegraphics[width=0.76\textwidth]{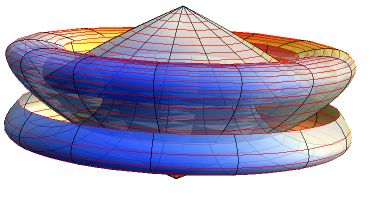}
  		\caption{}
  	\end{subfigure}
  	\begin{subfigure}[b]{0.45\linewidth}
  		\centering
  		\includegraphics[width=0.68\textwidth]{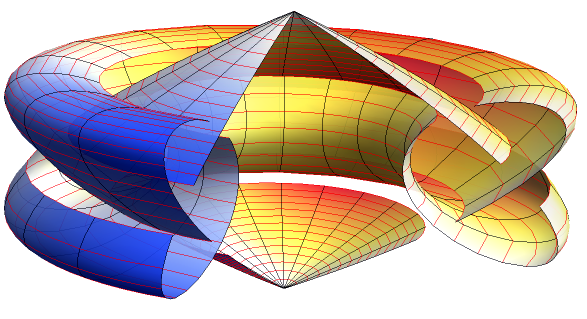}
  		\caption{}
  	\end{subfigure}  
  	\caption{A broken surface diagram of a spun trefoil knot}
  	\label{fig:brokendiagram}
  \end{figure}

To study a $2$-knot, we  project its image on a suitable $3$-dimensional space.  Let $F$ be a closed surface.
A map $g:F \rightarrow \R^{3}$ is called a {\it generic map} or a map in a {\it general position} if every point in $g(F)$ is either a regular point, a double point, a triple point or a branch point (Figure \ref{fig:singular points}).
\begin{definition}
	Choose a vector $v \in \R^{4}$ and let us take a projection of $\R^{4}$ onto some hyperplane $H$ orthogonal $v$ by a projection map $\pi:\R^{4} \rightarrow H$. A 2-knot $F$ is in general position with respect to $\pi$ or $\pi$ is a {\it generic projection} of $F$ if $\pi\vert_{F}:F \rightarrow H$ is a smooth immersion.
\end{definition}
\begin{remark} 
In general, a regular projection for a given $2$-knot may not exist. However, it has been proved \cite{Glu} that for every $2$-knot $F$ in $\R^4$, there exists a $2$-knot $\tilde{F}$ which is ambient isotopic to $F$ and for which a regular projection exists.
\end{remark}
%\begin{figure}[h]
%	\centering
%	\begin{subfigure}[c]{0.4\linewidth}
%		\centering
%		\includegraphics[width=0.8\textwidth]{doublearcs.png}\label{fig:doublearcs}
%		\caption{For a double point arc }
%	\end{subfigure}
%	\quad
%	\begin{subfigure}[c]{0.4\linewidth}
%		\centering 
%		\includegraphics[width=0.8\textwidth]{doublecircle.png}\label{fig:doublecircle}
%		\caption{For a double point circle}
%	\end{subfigure} 
%	
%	\caption{Neighborhood of double point sets }
%	\label{fig:doublecurves}
%\end{figure}
In a regular projection, the image of $S^2$  in $\R^3$ is an immersion. Hence, the derivative map is injective at each point of $S^2$. Thus, the singular set (points where the map fails to be one to one) will consist of a finite number of double curves and isolated triple and branch points.	
A broken surface diagram \cite{CarSai} of a $2$-knot consists of the image $\pi(F)$ under a generic projection map $\pi$ equipped with the over/under information on each double point set and on each triple point with respect to the direction of the vector $v$.
  For a singular point, we can distinguish the sheets (i.e. the 2-disks) by the relative height with respect to the direction of the vector $v$ as upper-lower (for double points and branch points), upper-middle-lower (for triple points). In case of branch point, two sheets converge at the branch point. A broken surface diagram is obtained by cutting the lower sheet consistently along all double point curves. The cutting follows naturally for the triple and the branch points as the triple points are contained in the double point set and the branch points are in the closure of the double points set. The neighborhood of the singular points after removing a strip from the lower sheet is shown in Figure \ref{fig:broken} and a broken surface diagram of a knotted sphere (spun trefoil knot) is shown in Figure \ref{fig:brokendiagram}.
  
  It is hard to visualize a 2-dimensional knot in the 4-dimensional space. So, we consider its projection to a 3-dimensional subspace.
  But it is sometimes hard to draw the whole surface diagram as some parts are behind
  other parts. To avoid that we can use movies or motion pictures where we cut the surface diagram by parallel hyperplanes. Then each intersection gives a 1-dimensional manifold.
  Let $F \subset \R^{4}$ be a given knotted surface and take a vector $v \in \R^{4}$. Let $\{\R^{3}_{t}\}_{t\in \R}$ be family of parallel hyperplanes, each of which is orthogonal to the vector $v$, containing the point $tv$. Now, $F_{t} = F \cap \R^{3}_{t}$ ia always a 1-dimenisonal manifold i.e a classical link or empty. Since, $F$ is compact, there exists an interval $[a,b]$ such that $F \subset \cup_{t\in[a,b]}\,\R^{3}_{t}$.

  Let $p_{v}: \R^{4} \rightarrow \R_{v}$ be the projection map onto the line $\R_{v}=\{tv\in \R^{4}\,\vert \,t\in \R\}$. Consider the restriction map $p_{v}\vert_{F}:F \rightarrow \R_{v}$.
  By Morse theory, non-degenrate critical points of $p_{v}\vert_{F}$ are isolated and the intersections at those points are of our interest as the knot type does not change at any point other than the critical points. 
  \begin{definition}
  	For the Morse function on $F$ defined above by restricting the projection of $\R^{4}$ onto $\R_{v}$, the one parameter family $\{F_{t} \subset \R^{3}_{t}\}_{t\in[a,b]}$ is called the {\it motion picture} or a {\it movie} of the surface link $F$. And each $F_{t}$ for a $t\in \R$ is called a {\it still} of the movie.
  \end{definition}
  \begin{remark}
  	Since the knot type changes at only at finitely many critical points, a finite number of stills are sufficient to describe the knotted surface. Therefore, although the whole famliy $\{F_{t} \subset \R^{3}_{t}\}_{t\in[a,b]}$ is called the movie, we can take only a finite set of stills as the motion picture or a movie of the surface link.
  \end{remark}
  A detailed exposition on $2$-knots can be found in  \cite{CarSai} and \cite{Kam}.
 \subsection{Construction of knotted spheres:}\label{sec:sphere} 
 In this section, we discuss the constructions to obtain knotted surfaces in 4-space such as spinning and twist spinning.
 \subsubsection{Spinning construction}\label{sec:spin}
\begin{figure}[h]
	\centering
	\begin{subfigure}[c]{0.4\linewidth}
		\centering
		\includegraphics[width=0.65\textwidth]{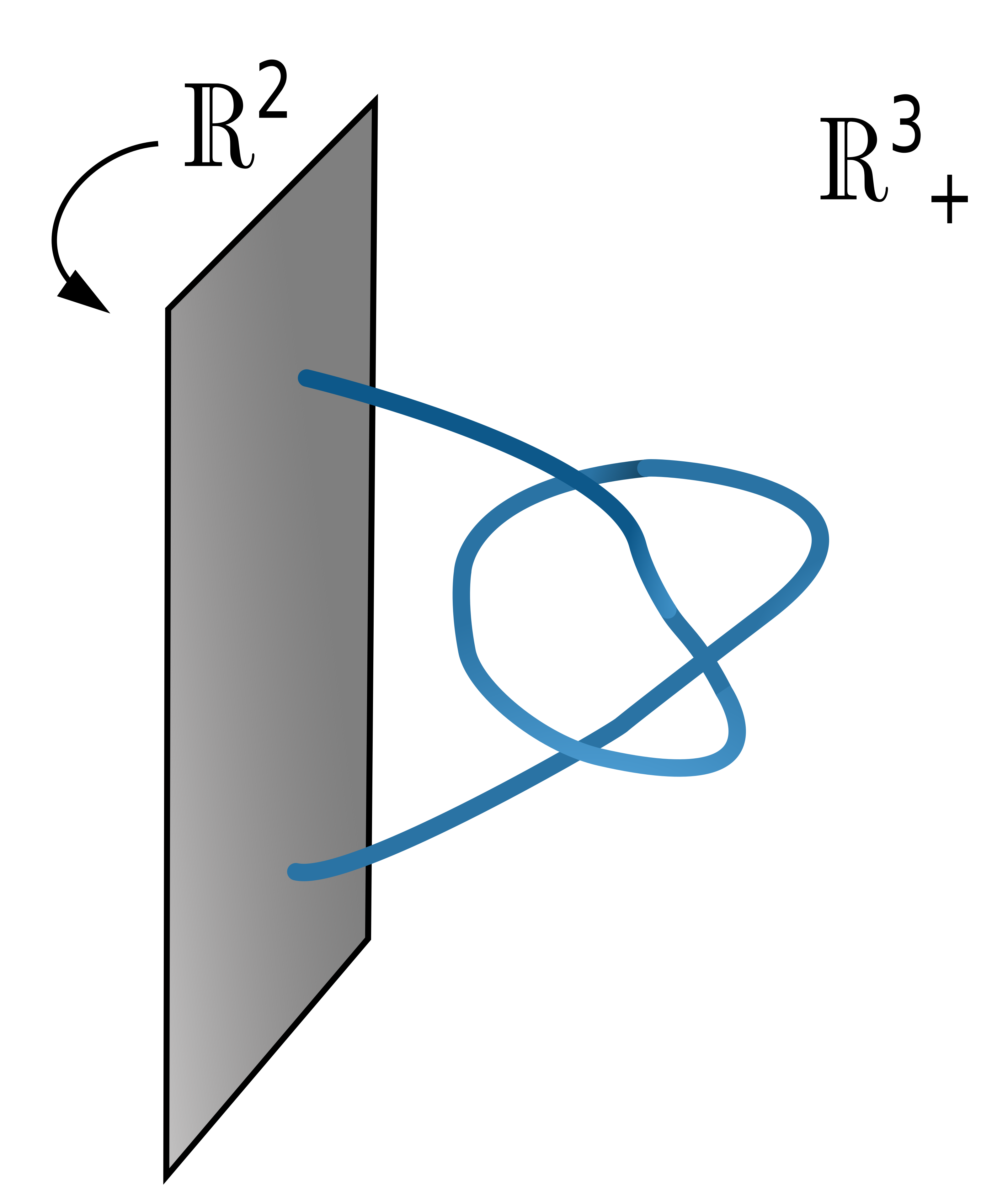}
		\caption{Spinning}
		\label{fig: arc in upper space}
	\end{subfigure}
	\begin{subfigure}[c]{0.4\linewidth}
		\centering
		\includegraphics[width=0.65\linewidth]{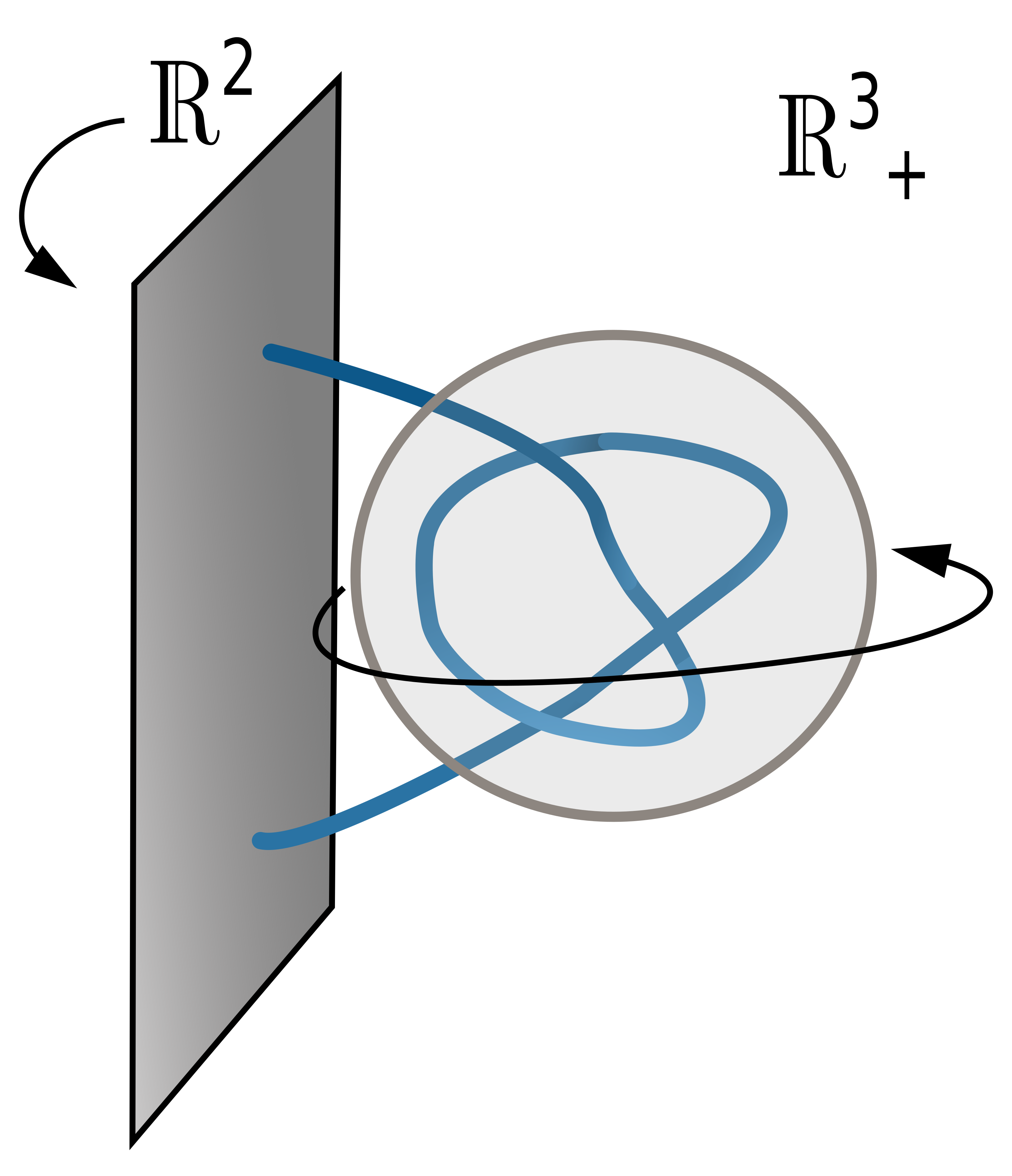}
		\caption{Twist spinning}
		\label{fig:twistspinning}
	\end{subfigure}
	\caption{Knotted sphere constructions}
\end{figure}

E. Artin \cite{Art} introduced a way to construct $2$-knots by spinning a knotted arc embedded in the half-space $\R^{3}_{+}$ around $\R^{2}$. The resulting $2$-knots are called {\it spun 2-knots}. Spun knots represent the simplest class of $2$-knots. Their construction is described as follows. In $\R^{4}$, consider the upper half space \[\R^{3}_{+}=\LC(x, y, z, 0): z\geq 0 \RC\] 
with the boundary
\[\partial \R^{3}_{+}=\R^{2}=\LC(x, y, 0, 0)\RC.\]
Now, locus of a point $x=(x, y, z, 0)$ in $\R^{3}_{+}$, rotated in 4-space about $\R^{2}$ is given by
\[\LC(x, y, z \,\cos \theta, z \,\sin \theta) \MV  0 \leq\theta\leq 2\pi \RC.\]
To construct a spun 2-knot, we choose a properly, locally flatly embedded arc $k$  in $\R^{3}_{+}$ i.e $k$ is embedded in $\R^{3}_{+}$ locally flatly and intersects $\partial \R^{3}_{+}$ transversely only at the endpoints (Figure \ref{fig: arc in upper space}). Then we spin $\R^{3}_{+}$ along $\R^{2}$ by $360^{o}$ and the continuous locus of the arc $k$ produces the spun 2-knot given by
\[\LC(x, y, z \,\cos \theta, z \,\sin \theta) \MV (x, y, z)\in k, 0 \leq\theta\leq 2\pi \RC.\]

\subsubsection{Twist spinning construction}\label{sec:twistspin}
E. C. Zeeman generalized Artin's spinning construction to {\it twist spinning} in 1965 \cite{Zee65}. In this case, we imagine the knotted part of $k$ inside a 3-ball as in Figure \ref{fig:twistspinning}. Now, to include twisting in the spinning construction we rotate the ball $d$ times around its own axis while rotating $\R^{3}_{+}$ around $\R^{2}$ once. Position of the knotted arc after $d$ twists should match its initial position after completing the rotation around $\R^{2}$. This way we get a 2-sphere in $\R^{4}$, called {\it d-twist spun 2-knot}. By definition, for $d=0$ we get a spun knot.

\subsection{Classical polynomial knots}\label{sec:polyknots}

It is already known that the ambient isotopy classes of classical knots ($S^1\subset S^3$) are in bijective correspondence with the ambient isotopy classes of long knots, i.e., the smooth embeddings of $\R$ in $\R^3$ which are proper and have asymptotic behaviour outside a closed interval. In 1994, the following theorems were proved regarding parameterization of long knots.
\begin{theorem} [\cite{Sha}, \cite{Shu}]
	
	For every long knot $K$ there exists a polynomial embedding $t\to (f(t),g(t), h(t))$ from $\R$ to $\R^3$  which is isotopic to $K$.
\end{theorem}

\begin{theorem}[\cite{Shu}]
	
	If two polynomial embeddings $\phi_0: \R\to \R^3$ and $\phi_1: \R\to \R^3$ represent isotopic knots then there exists a one-parameter family of polynomial embeddings $\phi_t: \R\to \R^3$ for each $t\in [0,1]$. We say that $\phi_0$ and $\phi_1$ are polynomially isotopic.
\end{theorem}

Both of these theorems stated above were proved using Weierstrass' Approximation \cite{Rud}. Later, concrete polynomial embeddings representing a few classes of knots were constructed, and their degrees were estimated (\cite{PraMis}, \cite{Mish09}). In all these constructions the main idea was to find a concrete embedding  to fix a suitable knot diagram and find real polynomials $f(t)$ and $g(t)$ such that the plane curve
$(x(t), y(t))= (f(t),g(t))$ provides the  projection of the diagram. For this curve, one can explicitly find the parameter pairs $(s_i,t_i)$ for which $f(s_i)=f(t_i)$ and $g(s_i)=g(t_i)$ for each $i$. Then, depending upon the over/under crossing information, the authors construct a polynomial $h(t)$ which satisfies $h(s_i)<h(t_i)$, at an undercrossing and $h(s_i)>h(t_i)$, at an overcrossing. Note that, in these embeddings, the knot is always projected
on the $XY$ plane, and the $Z$ co-ordinate provides the height function. All the crossings of the knot lie in the image of a closed interval $[a,b]$.

Here, we prove the following lemma, which will be useful in Section \ref{sec:sphere}.

\begin{figure}[h]
	\centering
	\includegraphics[width=0.35\textwidth]{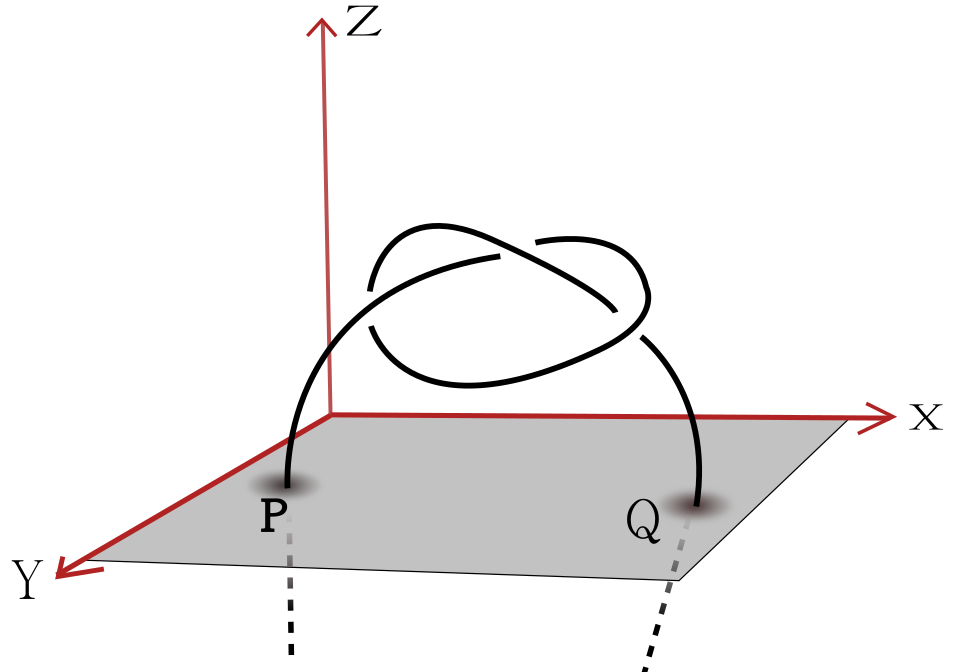}
	\caption{A long knot with knotting part in $\R^{3}_{+}$}
	\label{fig:long trefoil 1}
\end{figure}
\begin{lemma}\label{lem: knotted arc}
	Given a long knot  $K$, there exists a polynomial embedding  from $\R\to \R^3$ defined by $t\to (x(t),y(t),z(t))$  such that for some interval $[a,b]$,  $z(a)=z(b)=0$, $z(t)>0$ for $t\in (a,b)$ and there are no crossings outside $[a,b]$.
\end{lemma}
\begin{proof} We construct polynomials $f(t)$ and $g(t)$ such that the curve, defined by $(x,y)= (f(t),g(t))$ represents a projection of $K$ on the $XY$ plane. We find the parameters $(t_i,s_i)$ such that $f(s_i)=f(t_i)$ and $g(s_i)=g(t_i)$. We arrange $t_i$ and $s_i$ in ascending order. We choose $\mu_j$ between any two parameters which are giving rise to crossing points of the same kind (all over or all under) and consider $\tilde{h}(t)=\pi_j (t-\mu_j)$. This polynomial $\tilde{h}(t)$ provides the over/under crossing information. Since our crossings are finite, they lie within a closed interval say $[c,d]$. If the degree of $\tilde{h}(t)$ is odd, we can choose a real number $\mu$ outside $[c,d]$ in such a way that the polynomial $h_0(t)= -\tilde{h}(t) (t-\mu)$ still provides the over/under crossing information. Note that  $h_0(t)$ is an even degree polynomial. Now by adding a sufficiently large real number $R$, we can show that $h(t)=h_0(t)+R$ has exactly two real roots $a$ and $b$ outside $[c,d]$,  $h(t)>0$ for $t\in (a,b)$ and all the double points of the projection lie inside $[a,b]$ (Figure \ref{fig:long trefoil 1}). This completes the proof.
\end{proof}

\begin{figure}[h!]
	\centering
	\includegraphics[width=3.3cm]{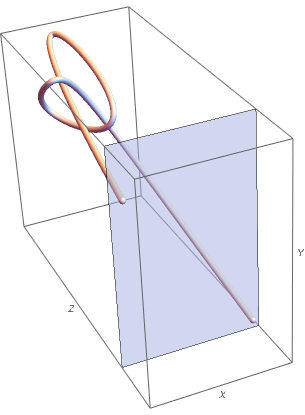}
	\caption{}
	\label{fig:one embedding}
\end{figure}

\begin{example}  
	As a demonstration, here is an embedding of the long trefoil knot given by
	\[t\to (t^3-3t, t^4-4t^2, -t^6+2t^5+4.24t^4-8.48t^3-3.24t^2+6.48t+12).\] For the interval $[a,b]=[-1.94871,2.40711]$, we get the knotted part in $\R^{3}_{+}$ and there are no crossings outside $[-1.94871,2.40711]$ (Figure \ref{fig:one embedding}).

\end{example}
\section{ Polynomial Representation of long $2$-knots }\label{sec:long}

In this section, we generalize the results in \cite{Shu} in the context of $2$-knots.

\begin{theorem}
	Every long $2$-knot given by a smooth embedding $\phi:\R^2 \rightarrow \R^4 $, has a polynomial representation.
\end{theorem}
\begin{proof}
	Let $F$ be long 2-knot in
	$\R^4$ given by a smooth embedding $\phi:\R^2 \rightarrow \R^4 $
	defined by $\phi(s,t)=\LP \alpha(s,t),\beta(s,t),\gamma(s,t),\delta(s,t)\RP$  such that $\Tilde{\phi} \equiv \LP \alpha,\beta,\gamma\RP:\R^2 \rightarrow \R^3$ is a generic immersion (having finite number of double point sets and triple points as the only singularities). 
	Up to equivalence, we can assume that the Jacobian matrix 
	\[J=
	\begin{pmatrix}
		\frac{\partial \alpha}{\partial s} & \frac{\partial \alpha }{\partial t}\\[2mm]
		\frac{\partial \beta }{\partial s} & \frac{\partial \beta }{\partial t}\\[2mm]
		\frac{\partial \gamma }{\partial s} & \frac{\partial \gamma }{\partial t}
	\end{pmatrix}
	\]
	has rank $2$ outside some closed region $\I_{1}=[a,b] \times [a,b]$.
	
	Since we have finitely many double point sets or triple points in the image of $\Tilde{\phi}$, we can choose another closed region $\I_{2}=[c,d] \times [c,d]$ such that $\Tilde{\phi}\LP \I_{2} \RP$ contains all the singularities. Let,
	$\M=[m_{1}, m_{2}] \times [m_{1}, m_{2}]$ such that $\M$ contains both the rectangles $ \I_{1}$
	and $\I_{2}$. Let $\phi(\M) $ be contained inside a $4$-ball of radius $r$ with 
	\[ \| \phi(m_{i}, m_{j}) \|= r \quad \forall \, i,j=1,2.\]
	Let $\mathcal{N}=[n_{1}, n_{2}] \times [n_{1}, n_{2}]$ be such that  $\phi(\mathcal{N}) $ is contained inside a $4$-ball of radius $2r$ with 
	\[ \| \phi(n_{i}, n_{j}) \|=2r \quad \forall \, i,j=1,2.\]
	By a smooth reparameterization, we can assume that 
	\[[m_{1}, m_{2}]=\LS -\frac{1}{2},\frac{1}{2}\RS \quad \text{and } \quad [n_{1}, n_{2}]=[-1,1].\]
	Thus, $\phi\LP \M \RP$ and $\phi(\mathcal{N})$ are contained inside balls of radius $r$ and $2r$ respectively with 
	\[\Big\| \phi\LP \pm\frac{1}{2},\pm\frac{1}{2}\RP \Big\|= r \quad
	\text{and} \quad
	\| \phi(\pm 1,\pm 1) \|= 2r\]
	and the Jacobian matrix has rank $2 $ outside $\M$.
	Now, consider the restriction of $\phi$ to $\I_{1}$ i.e 
	\[\phi\vert_{\I_{1}}:\I_{1} \rightarrow \R^4.\]
	Since the set of embeddings from a compact Hausdorff manifold to  any manifold forms an open set in the set of all smooth maps with $C^{1}$-topology, there exists an $\varepsilon_{0} > 0$ such that
	$\psi \in N(\phi,\varepsilon_{0})$ i.e.  $\psi$ is an embedding of $ \I_{1}$ in $\R^4$,
	where
	\[N(\phi,\varepsilon_{0})=\LC \psi: \underset{(s,t) \in \I_{1} }{\sup} \{ \|\psi(s,t)-\phi(s,t)\|,\|\psi'(s,t)-\phi'(s,t)\|\} < \varepsilon_{0}\RC.\]
	Let $\varepsilon < \min\{R/2,\varepsilon_{0}\}.$
	For this $\varepsilon$, we can take an $\frac{\varepsilon}{2}$-approximation $\psi$ using the Bernstein polynomial \cite{Sta}
	or the Weierstrass' Approximation 
	of two variables inside the square $\mathcal{N}$ \cite{Sta}. Let $\psi_{1} \equiv (x_{1},y_{1},z_{1},w_{1})$ be the approximation,where $x_{1},y_{1},z_{1}$ and $w_{1}$ are polynomials in two variables $s,t$.
	Then all the partial derivatives of the polynomials $x_1,y_1,z_1$ and $w_1$ are greater than $ (1-\frac{\varepsilon}{2})$ in the region $\mathcal{N} \setminus \M$.
	
	Now, for $\delta_{i}\in (0,\varepsilon/2),i=1,2$, we can choose $\N \in \mathbb{Z}^{+}$ large enough so that
	\begin{align*}
		\psi & \equiv  \Big(  x_{1}+\frac{\delta_{1}}{2\N+1}s^{2\N+1}+\frac{\delta_{2}}{2\N+1} t^{2\N+1},\,
		y_{1}+\frac{\delta_{1}}{2\N+1}s^{2\N+1}+\frac{\delta_{2}}{2\N+1}t^{2\N+1}, \\
		& \qquad  z_{1}+\frac{\delta_{1}}{2\N+1}s^{2\N+1}+\frac{\delta_{2}}{2\N+1}t^{2\N+1},\,
		w_{1}+\frac{\delta_{1}}{2\N+1}s^{2\N+1}+\frac{\delta_{2}}{2\N+1}t^{2\N+1} \Big)  \\
		& \equiv (x,y,z,w).
	\end{align*}
	Then $\psi$ is an $\varepsilon$-approximation of $\phi$ inside $\mathcal{N}$ such that the derivative of $\psi$  is positive outside the square  $\mathcal{N}$ as each partial derivative 
	\[\frac{\partial x}{\partial s},\frac{\partial x}{\partial t},\frac{\partial y}{\partial s},\frac{\partial y}{\partial t},\frac{\partial z}{\partial s},\frac{\partial z}{\partial t},\frac{\partial w}{\partial s},\frac{\partial w}{\partial t}\]
	are positive outside $\mathcal{N}$. Now, inside a compact region, this polynomial map $\psi$ approximates the knot. Also the presence of large odd degree  term with a very small coefficient ensures that the long $2$-knot becomes asymptotically flat outside the compact region.
	
	It remains to show that $\psi:\R^2 \rightarrow \R^4$ is an embedding. 
	$\phi:\R^2 \rightarrow \R^4$ being an embedding implies it is an injective immersion. That means the Jacobian $J(\phi;s,t)$ has rank $2$.  Now, in $\mathcal{N}$, \[\|\psi'(s,t)\| >\|\phi'(s,t)\|-\varepsilon/2\] which implies that the Jacobian $J(\psi;s,t)$ has rank $2$ in $\mathcal{N}$.
	And, outside 
	$\mathcal{N}$, all partial derivatives of $\psi$ are positive which implies that the Jacobian $J(\psi;s,t)$ has rank $2$ outside $\mathcal{N}$ as well. Hence, $\psi$ is an immersion. 
	Since the Jacobian matrix for $\psi$ has rank $2 $ both inside and outside $\mathcal{N}$, by inverse function  theorem it follows that $\psi$ is injective.
	This completes the proof.
\end{proof}
\begin{definition}
	Two polynomial embeddings $\psi_{0} \equiv (x_{0}, y_{0}, z_{0}, w_{0})$, $\psi_{1} \equiv (x_{1}, y_{1},$ $z_{1}, w_{1})$ of $\R^2$ in $\R^4$ are said to be \it{polynomially homotopic} or \it{P-homotopic} if there exists a homotopy  $F:\R^2 \times I \rightarrow \R^4$ between $\psi_{0}$ and $\psi_{1}$ such that $F_{u}:=F(s,t,u)$ is a polynomial map for each $u \in [0,1]$. In case each $F_{u}$ is a polynomial embedding, then $F$ is called a $P$-isotopy, and the embeddings $\psi_0$ and $\psi_1$ are said to be polynomially isotopic.
\end{definition}
\begin{notation}
	For $\varepsilon \in \R^{+},\N \in \mathbb{Z}^{+}$ and $\phi \equiv (x,y,z,w):\R^2 \rightarrow \R^4$, let us define 
	$\phi_{\varepsilon,N}:\R^2 \rightarrow \R^4$ given by \[\phi_{\varepsilon,N}(s,t)=\LP x(s,t), y(s,t), z(s,t)+\varepsilon \,s^{2\N+1}, w(s,t)+\varepsilon\,t^{2\N+1}\RP.\]
\end{notation}
We prove the following lemma.

\begin{lemma}\label{lem:odd deg}
	Let $\phi =(x,y,z,w):\R^2 \rightarrow \R^4$ be a polynomial embedding such that the map $(x,y,z):\R^2 \rightarrow \R^3$ is a generic immersion. Then for each $\N \in \mathbb{Z}^{+},$ there exists an $\varepsilon>0 $ such that $\phi$ and $\phi_{\varepsilon,N}$ are polynomially isotopic.
\end{lemma}

\begin{proof}
	To prove the lemma we show that the maps $F_{u}:\R^2 \rightarrow \R^4$ given by
	\[F_{u}(s,t)=(x(s,t), y(s,t), z(s,t)+u\varepsilon \, s^{2\N+1}, w(s,t)+u\varepsilon\,t^{2\N+1})\]
	are embeddings for each $u \in [0,1]$.
	
	The Jacobian for each $F_{u}$ 
	\[J=
	\begin{pmatrix}
		\dfrac{\partial x}{\partial s} & \dfrac{\partial x }{\partial t}\\[2mm]
		\dfrac{\partial y }{\partial s} & \dfrac{\partial y }{\partial t}\\[2mm]
		\dfrac{\partial z }{\partial s}+(2\N+1)u \varepsilon s^{2N} & \dfrac{\partial z }{\partial t}\\[2mm]
		\dfrac{\partial w }{\partial s} & \dfrac{\partial w }{\partial t}  +(2\N+1)u \varepsilon t^{2N}\\
	\end{pmatrix}
	\]
	has rank $2$ since the Jacobian for $\phi$ has rank $2$. Hence, $F_{u}$ is an immersion for each $u \in [0,1]$. Next, we show that $F_{u}$ is injective for each $u \in [0,1]$. Let us use the notation
	\[x(s_{i}, t_{i})=x_{i,i}.\]
	Now, consider the set 
	\[S=\LC (s_{1}, t_{1}),(s_{2}, t_{2})\in \R^2 \, \MV \,(x_{1,1},y_{1,1})= (x_{2,2},y_{2,2})\, \text{for}\, (s_{1}, t_{1}) \neq (s_{2}, t_{2})\RC.\]
	Since $\phi$ is an embedding, 
	\[z_{1,1} \neq z_{2,2} \quad \text{when} \quad w_{1,1}=w_{2,2} \]
	or
	\[w_{1,1} \neq w_{2,2} \quad \text{when} \quad z_{1,1}=z_{2,2}. \]
	We consider the first case. Now, we choose an $\varepsilon$ such that
	\[0 < \varepsilon < \underset{(s_{1}, t_{1}),(s_{2}, t_{2})\in S} {\text{min}} 
	\LC \frac{\Big\vert z_{1,1}-z_{2,2} \Big\vert}{\Big\vert s_{1}^{2\N+1}-s_{2}^{2\N+1} \Big\vert } \RC.\]
	Then, notice that this $\varepsilon$ satisfies
	\[z_{1,1}+u \varepsilon s_{1}^{2\N+1} \neq z_{2,2}+u \varepsilon s_{2}^{2\N+1}.\]
	for all $(s,t) \in S$ and $u\in [0,1]$. An analogous argument works for the second case.
	Thus, $F_{u}$ is injective for all $u \in [0,1]$. This completes the proof.
\end{proof}
This lemma allows us to choose a polynomial embedding with a large odd degree in third and fourth coordinates that keeps the Jacobian positive outside a compact region, say $
\M$.

\begin{theorem}
	Let $\psi_{0} \equiv (x_{0},y_{0},z_{0},w_{0})$ and $\psi_{1} \equiv (x_{1},y_{1},z_{1},w_{1})$ be two polynomial embeddings of $\R^2$ in $\R^4$ which represent the isotopic $2$-knots. Then $\psi_{0}$ and $\psi_{1}$ are polynomially isotopic.
\end{theorem}
\begin{proof}
	We can always choose $\psi_0$ and $\psi_1$ such that $\Tilde{\psi_{0}}=(x_0,y_0,z_0):\R^2 \rightarrow \R^3$ and $\Tilde{\psi_{1}}=(x_1,y_1,z_1):\R^2 \rightarrow \R^3$ are generic immersions and both $\psi_0$ and $\psi_1$ are asymptotic outside a compact region say $\M$.
	Now, let $\A$ be the common compact region inside which both $\Tilde{\psi_{0}}$ and $\Tilde{\psi_{1}}$ have their double curve and triple points.
	Choose $\mathcal{N} \supset \M \cup \A$ and image of $\mathcal{N}$ under $\psi_{0}$ and $\psi_{1}$
	are $4$-balls of radius $r_{1}$ and $r_{2}$ respectively.
	Let $r=\text{max}\{r_{1},r_{2}\}$. We can choose $\mathcal{N}=[-\frac{1}{2},\frac{1}{2}]\times [-\frac{1}{2},\frac{1}{2}]$.
	Since $\psi_{0}$ and $\psi_{1}$ represent the same knot type, there exists an isotopy $F:\R^2 \times I \rightarrow \R^4$ such that $F(s,t,0)=\psi_{0}(s,t)$ and $F(s,t,1)=\psi_{1}(s,t)$. Now consider the restriction of the isotopy map $F$ on $\I_{1}\times I$ where $\I_{1}=[-1,1]$. Each map $F_{u}:\I_{1}\rightarrow \R^4$ is on compact region $\I_{1}$ where we can choose a Weierstrass approximation of this map using Bernstein polynomial with two variables inside an $\varepsilon_{0}$ neighbourhood of $F_{u}$ for some $\varepsilon_{0} >0$ and adding higher odd order terms to the co-ordinates we can get polynomial embeddings $P_{u}:\I_{1}\rightarrow \R^4$ for each $u \in [0,1]$ whose derivative is positive outside. Therefore, $P_{0}$ and $P_{1}$ are polynomially isotopic. Now we define isotopies between $\psi_{0}$ and $P_{0}$ given by $H:\R^2 \times I \rightarrow \R^4$ where
	\[H(s,t,u)=(1-u)\psi_{0}(s,t) + u P_{0}(s,t).\]
	
	Now each $H_{u}$, defined by $H_{u}(s,t)=(1-u) \psi_{0}(s,t)+u P_{0}(s,t)$, is an $\varepsilon$- approximation of $\psi_{0}$ inside $\I_{1}$ for our choice of $\varepsilon$, $P$ and $\| H_{u} \|$ is increasing outside $[-\frac{1}{2},\frac{1}{2}]\times [-\frac{1}{2},\frac{1}{2}] $. We can use  similar arguments as in Theorem $6.0.1$ to prove that $H_{u}$ is an embedding of $\R^2$ in $\R^4$. Thus, we can define a polynomial isotopy between $\psi_{0}$ and $P_{0}$. Similarly we can show that $\psi_{1}$ and $P_{1}$ are polynomially isotopic. Hence $\psi_{0}$ and $\psi_{1}$ are polynomially isotopic.
	
\end{proof}

\section{Polynomial Representation of some families of $2$-knots:}\label{ch:sphere}

In Section \ref{sec:long}, we showed that every long $2$ knot can be realized as an image of $\R^2$ in $\R^4$ under a polynomial map, which is an embedding. By one-point compactification we get a compact $2$ knot, i.e.,
a knotted $S^2$ in $S^4$ and the map will no longer remain polynomial. However, in this section, we show that the spun knots and the twist spun knots can be realized as a compact image of a polynomial map.
In this section, we describe a method to write polynomial parameterizations spun knot and twist spun knot constructions, denoted by $spun(K)$ and the $d$-$twist spun (K)$, respectively. Section \ref{sec:sphere} shows how these knotted spheres are constructed with the help of certain rotation operations performed on a knotted arc. section we find functions representing the rotation operations involved in these constructions. Both of these knotted spheres require a locally flat, properly embedded knotted arc $k$ in ${\R^3}_{+}$ with endpoints on the boundary $\R^{2}$. By Lemma \ref{lem: knotted arc}, we can choose a polynomial parameterization of a long knot given by 
\[\phi(t) \rightarrow \LP f(t), g(t), h(t) \RP \]
such that $h(a)=0=h(b)$ and $h(t) >0$ for all $t \in (a,b)$. Therefore, $\phi([a,b])$ represents the required knotted arc $k$. We use this parameterization for the knotted arc $k$. 

The rotation operations applied on $k$ result in embeddings represented by smooth functions, not all of which are polynomials. Therefore, to get a polynomial parameterization, we use the Chebyshev approximations \cite{Abu} of these smooth functions inside a suitable compact domain. Section \ref{mathematica}, a {\it Mathematica} program is provided that is used to compute the Chebyshev approximations of $\cos dx$ and $\sin dx$ in any closed interval where $d \in \mathbb{N}$. For example,
the Chebyshev approximations $C(x)$ and $S(x)$ of $\cos x$ and $\sin x$ in $[0,2\pi]$ respectively are as follows:
\begin{align*}
	C(x) &= \,-0.0000193235 x^8+0.000485652 x^7-0.00399024 x^6+0.0081095 x^5 \\
	&+0.0265068 x^4 +0.0163844 x^3-0.509175 x^2+0.00205416 x+0.999921, \\[0.5mm]
	S(x) &= \,8.73651067430188 \times 10^{-19} x^8+0.000144829 x^7-0.00318496 x^6\\
	&+0.0220637 x^5-0.0322337 x^4-0.125592 x^3-0.0257364 x^2\\
	&+1.00614 x-0.000238495.
\end{align*}
For better clarity, we present the parameterizations of the spun knots and the $d$-twist spun knots in separate sections.

\subsection{Polynomial parameterization of spun $2$-knot}\label{sec:spunPoly}

\begin{theorem}\label{th:spun}
	Given a classical knot $K$, there exist polynomials $f(s,t)$, 
	$g(s,t), h(s,t)$ and $ p(s,t)$ in two variables 
	$t$ and $s$ such that for some interval $[a,b]$ the image of  $[a,b]\times [0,2\pi]$ under the map $\phi: \R^2 \to \R^4$ defined by $$\phi ((s,t))= \LP f(s,t), g(s,t), h(s,t), p(s,t) \RP$$ is isotopic to the spun of $K$.
\end{theorem}
\begin{proof}
	For the classical knot $K$, up to isotopy, let us choose a polynomial representation $\phi$ of $K$ such that  $\phi: [a,b] \rightarrow \R^{3}_{+}$ of the knotted arc $k$ is given by 
	\[\phi(t) \rightarrow \LP f(t), g(t), h(t)  \RP\] 
	such that $h(a)=0=h(b)$ and $h(t)>0$ for $t \in (a,b)$. Now, the spinning construction gives a map $F:[a,b]\times [0, 2\pi]\to \R^4$ defined by
	\[F(s,t)=\LP f(t),\, g(t),\,\, h(t)\,\cos s,\, h(t)\,\sin s \RP\] 
	that represents $spun(K)$.
	Using the point set topology argument, one can prove that the image of $F$ in $\R^4$ is homeomorphic to $S^2$. Now, we replace the trigonometric functions $\cos s$ and $\sin s$ with their Chebyshev approximations inside the interval $[0, 2\pi]$. Let the Chebyshev approximations of $\cos s$ and $\sin s$ inside the interval $[0,2\pi]$ be denoted by $C(s)$ and $S(s)$, respectively. Then by choosing $f(s,t)=f(t)$, $g(s,t)=g(t)$, $h(s,t)=h(t)\,C(s)$ and $p(s,t)=h(t)\,S(s)$, we obtain a polynomial parameterization of $spun(K)$. This completes the proof.
\end{proof}
\begin{figure}[h]
	\centering
	\includegraphics[width=0.65\textwidth]{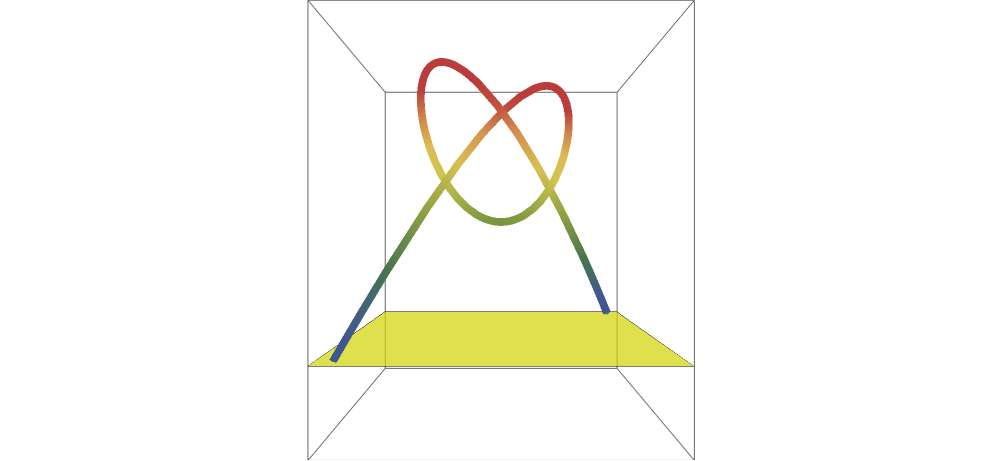}
	\caption{Knotted  arc of the long trefoil knot}
	\label{fig: long trefoil knot}
\end{figure}
\begin{example}[The spun trefoil knot] Let us take the following polynomial representation of the long trefoil knot given by
	\[ \LC \LP t^3-3t, t^5-10t , t^4-4t^2  \RP \MV t \in \R  \RC .\] 
	We change the function $(t^4-4t^2)$ to $(-t^4+4t^2+3)$
	whose real roots are $t= -2.1554 $ and $t=2.1554 $. Therefore, the knotted arc $k$ (Figure \ref{fig: long trefoil knot}) is given by
	$$ \LC \LP t^3-3t, t^5-10t , -t^4+4t^2 +3 \RP \MV t \in [-2.1554,2.1544] \RC .$$ 
	The spun trefoil knot is given by 
	\[ \LC  \LP t^3-3t, t^5-10t , (-t^4+4t^2 +3)\,\cos s, (-t^4+4t^2 +3)\,\sin s \RP \RC\]
	for $t \in [-2.1554,2.1544], s\in [0, 2\pi]$. After replacing cosine and sine functions with $C(s)$ and $S(s)$ respectively, we get the final polynomial parameterization as follows
	\[\LC \LP t^3-3t, t^5-10t , (-t^4+4t^2 +3)\,C(s), (-t^4+4t^2 +3)\,S(s) \RP \RC\]
	for $t \in [-2.1554,2.1544], s\in [0, 2\pi]$. {\it Mathematica} plots for a projection of the spun trefoil are shown in Figure \ref{fig:Spun trefoil projection}.
	In Figure \ref{fig:a hyperplane section}, an inside view is also given where it can be seen that the knotted arc does not deform while spinning around $\R^{2}$.
	\begin{figure}[h]
		\centering
		\begin{subfigure}[c]{0.4\linewidth}
			\centering
			\includegraphics[width=1.2\textwidth]{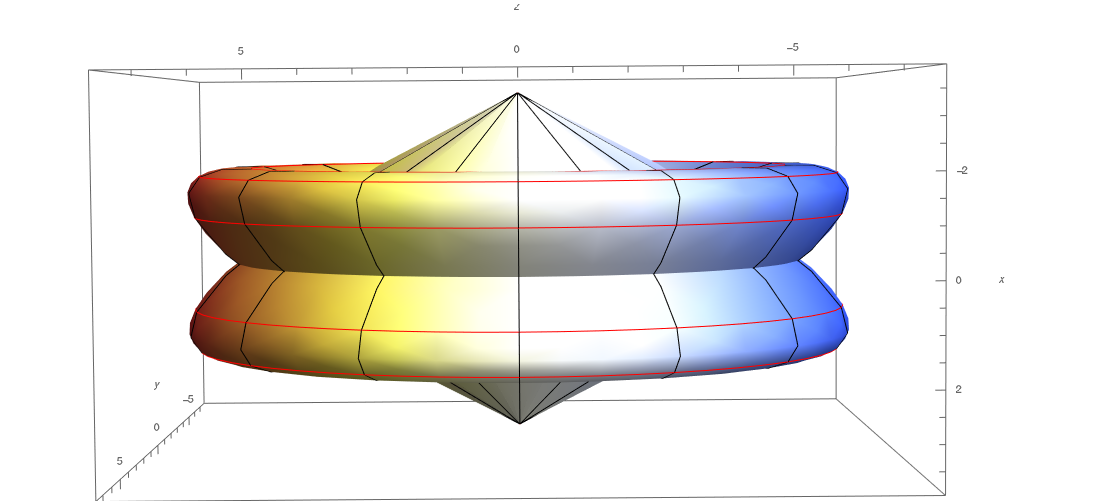}
			\caption{Projection on a hyperplane}
			\label{fig:Spun trefoil projection}
		\end{subfigure}
		\hspace{1.5cm}
		\begin{subfigure}[c]{0.4\linewidth}
			\centering
			\includegraphics[width=\textwidth]{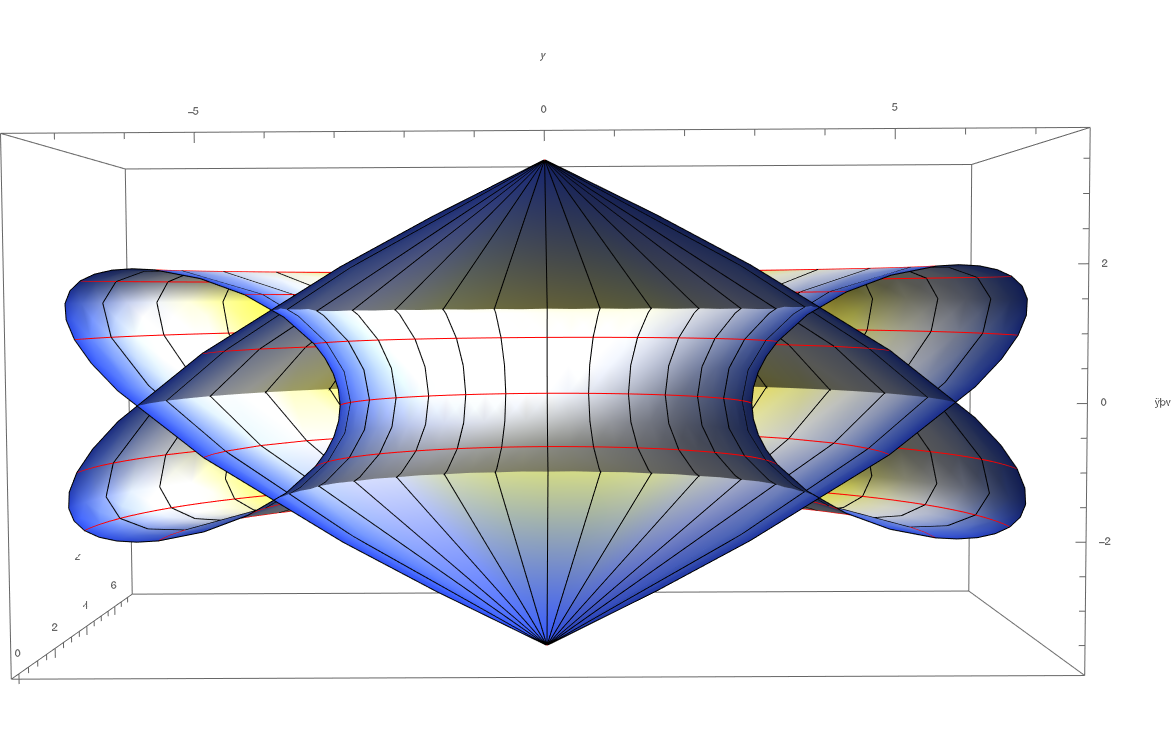}
			\caption{An inside view }
			\label{fig:a hyperplane section}
		\end{subfigure}
		\caption{The spun trefoil knot}
		\label{fig:spun trefoil}
	\end{figure}
\end{example}

\begin{figure}[h]
	\centering
	\includegraphics[width=0.38\textwidth]{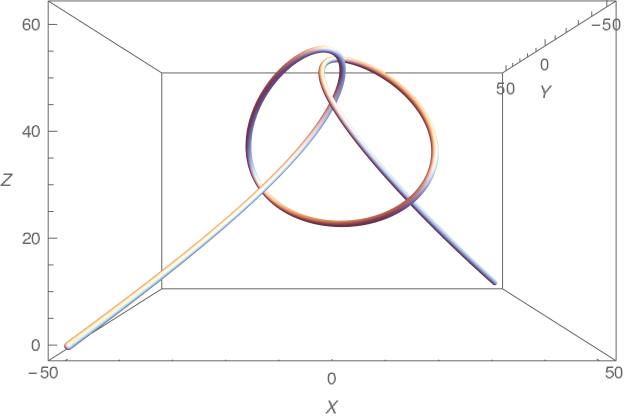}
	\caption{Knotted  arc of long figure eight knot}
	\label{fig: long figure eight knot}
\end{figure}

\begin{figure}[h]
		\centering
		\begin{subfigure}[b]{0.4\linewidth}
			\includegraphics[width=0.88\textwidth]{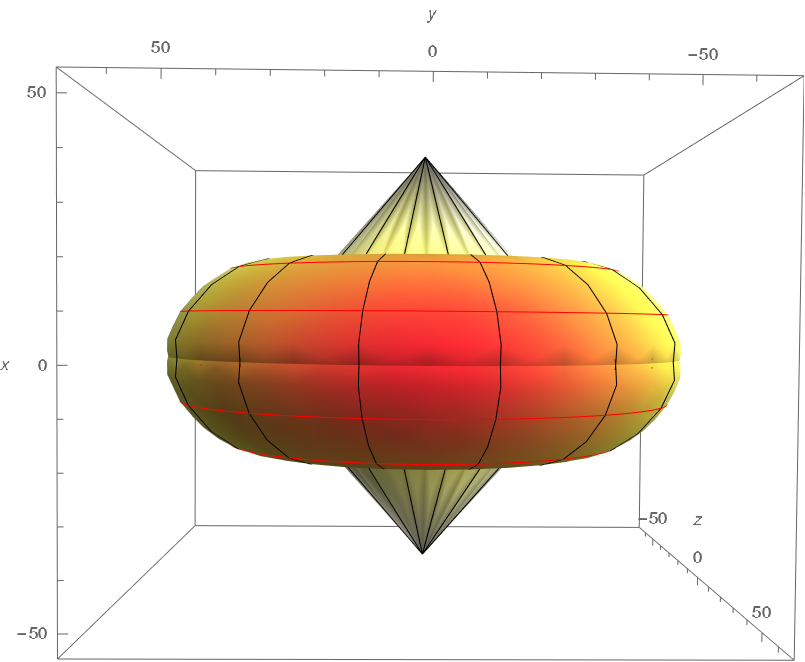}
			\caption{Projection on a hyperplane}
		\end{subfigure}
		\quad
		\begin{subfigure}[b]{0.4\linewidth}
			\includegraphics[width=0.9\textwidth]{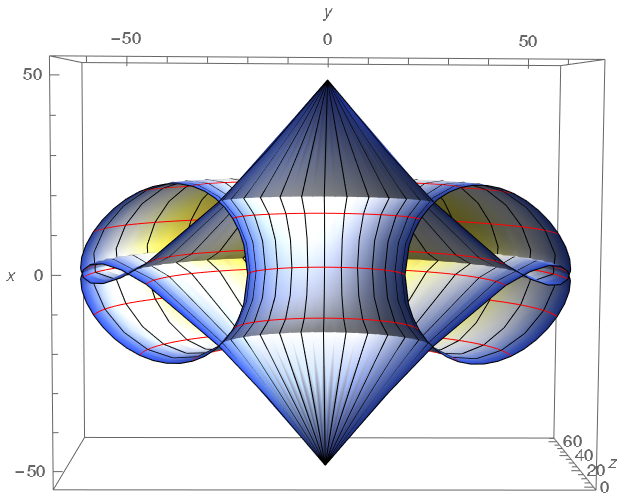}
			\caption{Inside view of the projection}
		\end{subfigure}
		\caption{The spun figure eight knot}
		\label{fig: spun figure eight}
	\end{figure}

\begin{example}[The spun figure eight knot]
	Using the parameterization for the long figure eight knot given in \cite{ANB} (Figure \ref{fig: long figure eight knot}) and following the same procedure as above, we obtain a polynomial parameterization of the spun figure eight knot given
	by $$\LC \LP f(t), g(t), h(t) \, C(s), h(t)\, S(s)\, \RP \MV t\in [-3.7934, 3.7934], s\in [0,2\pi]\RC$$
	where,
	\begin{align*}
		f(t) &=\frac{2}{5}  \LP t^2-7 \RP    \LP t^2-10 \RP   t, \\
		g(t) &= \frac{1}{10} t  \LP t^2-4 \RP    \LP t^2-9 \RP    \LP t^2-12 \RP  , \\
		h(t) &=(20-13 t^2-t^4).
	\end{align*}
	and $C(s)$, $S(s)$ are the Chebyshev approximations of $\cos s$ and $\sin s$ respectively. Figure \ref{fig: spun figure eight} shows the {\it Mathematica} plot of a projection of the spun figure eight knot on \textit{XZW-plane}. The {\it Mathematica} program is provided in Section \ref{spunfigure}.
	
\end{example}

%\begin{figure}[h]
%	\centering
%	\begin{
	%	\subfigure[long figure eight knot]{\includegraphics[width=5cm]{figure8projecton1.png}} 
	%	\subfigure[Spun figure eight knot]{\includegraphics[width=5cm]{fig8spun.png}} 
	%	\caption{}
	%	\label{fig:long figure eight and spun figure eight}
	%\end{figure}
	
	\subsection{Polynomial parameterization of $d$ twist spun knot}\label{th:twistPoly}
	\begin{theorem}
		Given a classical knot $K$, there exist polynomials $f(t,\theta)$, 
		$g(t,\theta), h(t,\theta)$ and $ p(t,\theta)$ in two variables 
		$t$ and $\theta$ such that for some interval $[a,b]$, the image of  $[a,b]\times [0,2\pi]$ under the map $\phi: \R^2 \to \R^4$ defined by $$\phi (t,\theta)= \LP \bar{f}(t,\theta), \bar{g}(t,\theta), \bar{h}(t,\theta), \bar{p}(t,\theta)\RP $$ is isotopic to the d-twist spun of $K$. 
	\end{theorem}
	\begin{proof}
		We start with the parameterization $\phi: [a,b] \rightarrow \R^{3}_{+}$ of the knotted arc $k$, given by 
		\[\phi(t) \rightarrow \LP f(t), g(t), h(t)  \RP\] 
		where $h(a)=0=h(b)$ and $h(t)>0$ for $t \in (a,b)$.
		The endpoints of $k$ on the boundary $\R^{2}$ is given by $(f(a),g(a),0)$ and $(f(b),g(b),0)$. Now we find an interval $[a',b'] \subset [a,b]$ where all the crossings of $K$ lie.
		
		\begin{figure}[h]
			\centering
			\includegraphics[width=0.3\textwidth]{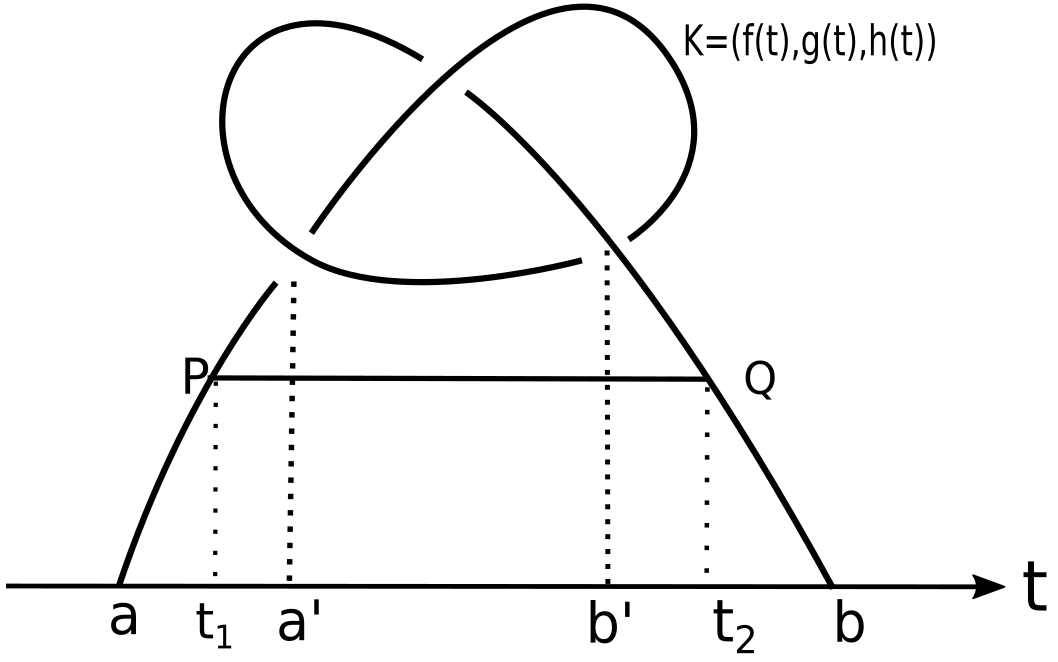}
			\caption{Rotating about $PQ$}
			\label{Rotating about $PQ$}
		\end{figure}

		In the twist spinning construction, the knotted part of the arc is contained in a 3-ball in $\R^{3}_{+}$ and is rotated around the axis of the ball. To achieve this first we choose the axis of rotation as a line segment $PQ$ parallel to the $XY$ plane joining two points on the knotted arc, say $P:=(f(t_{1}),g(t_{1}),h(t_{1}))$ and $Q:=(f(t_{2}),g(t_{2}),h(t_{2}))$, where $[a',b'] \subset [t_{1},t_{2}]$ and $h(t_{1})=h(t_{2})=c$ (Figure \ref{Rotating about $PQ$}). The value of $c$ is chosen so that the knotted part of the arc does not intersect the $XY$ plane while rotating around $PQ$. By Rodrigues' rotation formula \cite{Rod}, the rotation around $PQ$ is represented by the following matrix.
		\[\mathbf{R}'=\mathbf{T_{c}*R*T_{c}^{-1}}=\mathbf{T_{c}*R*T_{-c}},\]
		where $\mathbf{T_{c}}$ is the translation matrix along Z axis which will send $(x,y,z)$ to $(x,y,z+c)$ and $\mathbf{R}$ gives the rotation around the line $P'Q'$ on $XY$ plane, parallel to $PQ$, joining the points $P'=(f(t_{1}),g(t_{1}),0)$ and $Q'=(f(t_{2}),g(t_{2}),0)$ on the knotted arc.
		
		The rotation matrix $\mathbf{R}$ is defined as follows:
		If $\mathbf{v}$ is a vector in $\R^{3}$ and $\mathbf{k}$ is a unit vector describing an axis of rotation around which $\mathbf{v}$ rotates by an angle $\phi$ according to the right-hand rule, the Rodrigues formula for the rotated vector $\mathbf{v_{rot}}$ is given by
		\[\mathbf{v}_{rot}=\mathbf{v}\,\cos\phi +(\mathbf{k} \times \mathbf{v})\sin\phi+\mathbf{k}\cdot(\mathbf{k}\cdot \mathbf{v})(1-\cos\phi).\]
		In this case, $\mathbf{k}$ is along the line joining $(f(t_{1}),g(t_{1}),0)$ and $(f(t_{2}),g(t_{2}),0)$.
		So,\[\mathbf{k}=\Big(\frac{f(t_{2})-f(t_{1})}{N},\frac{g(t_{2})-g(t_{1})}{N},0\Big),\]
		where $N=\sqrt{(f(t_{2})-f(t_{1}))^{2}+(g(t_{2})-g(t_{1}))^{2}}.$
		
		For simplification, let us denote 
		\[f_{21}:=f(t_{2})-f(t_{1})\] 
		\[g_{21}:=g(t_{2})-g(t_{1}).\]
		Then\[\mathbf{k}=\Big(\frac{f_{21}}{N},\frac{g_{21}}{N},0\Big),\]
		where $N=\sqrt{f_{21}^{2}+g_{21}^{2}}.$
		
		Then the rotation matrix through an angle $\phi$ counterclockwise around the axis $\mathbf{k}$ is
		\[\mathbf{R=I+\sin\phi \; K+(1-\cos\phi)\;K^{2}},\]
		where
		\[ \mathbf{K}=
		\begin{pmatrix}
			0 & -k_{z} & k_{y} \\
			k_{z} &0 & -k_{x} \\
			-k_{y} & k_{x} & 0
		\end{pmatrix}
		=\begin{pmatrix}
			0 & 0 &\dfrac{g_{21}}{N} \\[0.5cm]
			0 & 0 & -\dfrac{f_{21}}{N} \\[0.5cm]
			-\dfrac{g_{21}}{N} & \dfrac{f_{21}}{N} & 0
		\end{pmatrix}.\]
		and the rotation matrix around $PQ$ is given by,
		\[\mathbf{R}'=\mathbf{T_{c}*R*T_{c}^{-1}}=\mathbf{T_{c}*R*T_{-c}},\]
		where
		\[\mathbf{R}=
		\begin{pmatrix}
			\dfrac{f_{21}^{2}\;+\;g_{21}^{2}\;\cos\phi}{N^{2}} & \dfrac{f_{21}\;g_{21}\;(1-\cos\phi)}{N^{2}} & \dfrac{g_{21}\;\sin\phi}{N} & 0\\[0.5cm]
			\dfrac{f_{21}\;g_{21}\;(1-\cos\phi)}{N^{2}} &   \dfrac{f_{21}^{2}\;\cos\phi\;+\;g_{21}^{2}}{N^{2}}  & -\dfrac{f_{21}\;\sin\phi}{N}& 0\\[0.5cm]
			-\dfrac{g_{21}\;\sin\phi}{N} & \dfrac{f_{21}\;\sin\phi}{N}&\cos\phi & 0 \\[0.5cm]
			0 & 0 & 0 & 1
		\end{pmatrix}\]
		and $\mathbf{T_{c}}$ is the translation matrix along Z-axis which sends $(x,y,z)$ to $(x,y,z+c)$, given by \[T_{c}=
		\begin{pmatrix}
			1 & 0 & 0 & 0 \\
			0& 1 & 0 & 0\\
			0 & 0 & 1 & c \\
			0 & 0 & 0 & 1
		\end{pmatrix}.\]
		Then,
		\[\mathbf{R}'=
		\begin{pmatrix}
			\dfrac{f_{21}^{2}\;+\;g_{21}^{2}\;\cos\phi}{N^{2}} & \dfrac{f_{21}\;g_{21}\;(1-\cos\phi)}{N^{2}} & \dfrac{g_{21}\;\sin\phi}{N} & -\dfrac{c\;g_{21}\;\sin\phi}{N}\\[0.5cm]
			\dfrac{f_{21}\;g_{21}\;(1-\cos\phi)}{N^{2}} &   \dfrac{f_{21}^{2}\;\cos\phi\;+\;g_{21}^{2}}{N^{2}}  & -\dfrac{f_{21}\;\sin\phi}{N} & \dfrac{c\;f_{21}\;\sin\phi}{N}\\[0.5cm]
			-\dfrac{g_{21}\;\sin\phi}{N} & \dfrac{f_{21}\;\sin\phi}{N}&\cos\phi& -c\;\cos\phi+c \\[0.5cm]
			0 &0 &0 & 1 
		\end{pmatrix}.\]
		Thus, the rotation around $PQ$ is given by the parameterization :
		\begin{equation}
			(t,\phi) \longrightarrow  \LC \LP f'(t,\phi),g'(t,\phi),h'(t,\phi) \RP\,\middle|\,
			\begin{aligned}
				a \leq &t \leq b \\
				0 \leq &\phi < 2\pi
			\end{aligned}\, \RC 
		\end{equation}
		where
		{\small \begin{align}
				f'(t,\phi)&=\frac{\LP f_{21}^{2}\;+\;g_{21}^{2}\;\cos\phi \RP\,\mathbf{f(t)}+f_{21}\;g_{21}\;(1-\cos\phi)\,\mathbf{g(t)}+N\;g_{21}\,\sin\phi\,(\mathbf{h(t)}-c)}{N^2}\\
				g'(t,\phi)&=\frac{f_{21}\;g_{21}\;(1-\cos\phi)\,\mathbf{f(t)}+\LP f_{21}^{2}\;\cos\phi\;+\;g_{21}^{2} \RP\,\mathbf{g(t)}+N\;g_{21}\,\sin\phi\,(c-\mathbf{h(t)})}{N^{2}}\\
				h'(t,\phi)&=\frac{-g_{21}\;\sin\phi\;\mathbf{f(t)}+f_{21}\sin\phi\;\mathbf{g(t)}+N
					\,\cos\phi\,\mathbf{h(t)}+N\,c\,(1-\cos\phi)}{N}
		\end{align} }
		and $N=\sqrt{f_{21}^{2}+g_{21}^{2}}.$
		
		Now, in the twist-spinning construction, there is no rotation outside the $3$-ball that contains the knotted part. This is same as not rotating the knotted arc outside the endpoints of $PQ$. We achieve this by combining the rotation functions with a bump function.
		\begin{figure}[h]
			\centering
			\includegraphics[width=0.3\textwidth,height=0.2\textwidth]{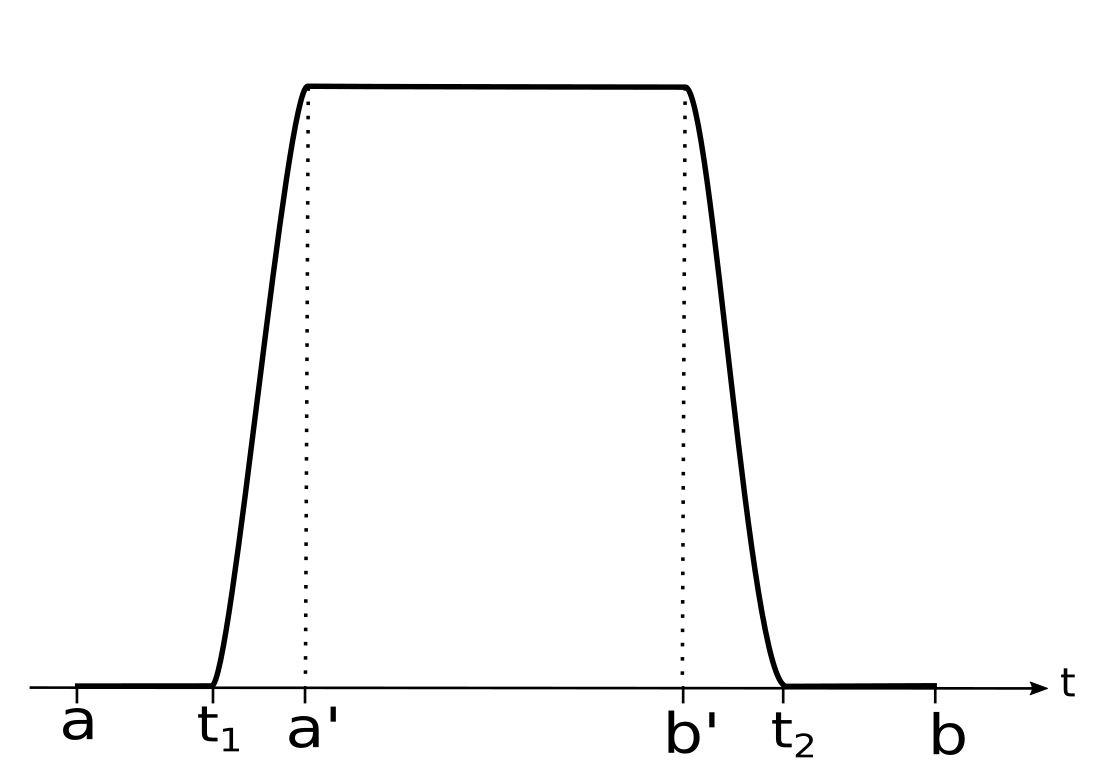}
			\caption{Bump function B(t). }
			\label{Bump function}
		\end{figure}
		
		We have, $ [a,b]= [a,t_{1}] \cup [t_{1},a'] \cup [a',b'] \cup [b',t_{2}] \cup [t_{2},b]$ (Figure \ref{Rotating about $PQ$}), where
		\begin{enumerate}
			\item $[a',b']$ contains all the crossings of the knot.
			\item $t_{1},t_{2}$ corresponds to the endpoints of $P$ and $Q$. So, $[a,t_{1}]$ and $[t_{2},b]$ correspond to the parts of the arc outside the axis $PQ$. 
		\end{enumerate}
		
		Therefore, we define a bump function which takes the value $1$ in $[a',b']$ and $0$ outside $[t_{1}, t_{2}]$. First, consider the function 
		\[F(t)=\begin{cases}
			e^{-\frac{1}{t}}\,, & t > 0 \\
			0 \,, & t \leq 0.
		\end{cases}\]
		Then define the bump function, $B(t)$ as follows:
		\[B(t)=\frac{F(d_{1}-t^{2})}{F(t^{2}-d_{2})+F(d_{1}-t^{2})},\]
		where $d_{1},d_{2} \in \R^{+} $ are chosen in such a way that 
		\[B(t)=\begin{cases}
			1\,, & t \in [a',b'] \\
			0\,, & t \in [a,t_{1}] \cup [t_{2},b]\\
			\in (0,1)\,,& \text{otherwise}.
		\end{cases}.\] 
		See Figure \ref{Bump function}.	Then, for $t \in [a,b]$, we define
		\begin{align*}
			\Tilde{f}(t,\phi)=& \,B(t)\,f'(t,\phi)+(1-B(t))\,f(t),\\
			\Tilde{g}(t,\phi)=& \,B(t)\,g'(t,\phi)+(1-B(t))\,g(t),\\
			\Tilde{h}(t,\phi)=& \,B(t)\,h'(t,\phi)+(1-B(t))\,h(t) .
		\end{align*}
		
		Now, the knotted arc is rotating $d$ times around the $PQ$ axis in $\R^{3}_{+}$ while rotating around the $XY$ plane in $\R^{4}$. Therefore, the rotation angle around $PQ$ is $d$ times the rotation angle around the $XY$ plane i.e. $\phi=d\,\theta.$
		
		Hence, the parameterization for d-twist spun knot is given by 
		\begin{equation}
			(t,\theta) \longrightarrow  \LC \LP \Tilde{f}(t,d\theta)\,,\,\Tilde{g}(t,d\theta)\,,\,\Tilde{h}(t,d\theta)\,\cos\,\theta\,,\,\Tilde{h}(t,d\theta)\,\sin\,\theta  \RP \MV
			\begin{aligned}
				a \leq &t \leq b \\
				0 \leq &\theta < 2\pi
			\end{aligned}  \RC .
		\end{equation}

		Here, $\Tilde{f},\Tilde{g},\Tilde{h}$ are linear combinations of the functions $f,g,h$, cosine, sine and the bump function $B(t)$. Replacing cosine, sine and $B(t)$ with their corresponding Chebyshev polynomial approximations in $[0,2\pi]$, we obtain a polynomial parameterization of d-twist spun knots given by 
		\begin{equation}
			\LC \LP \bar{f}(t,\theta)\,,\,\bar{g}(t,\theta)\,,\,\bar{h}(t,\theta)\,,\,\bar{p}(t,\theta)  \RP \MV a \leq t \leq b, \, 0 \leq \theta < 2\pi
			\RC .
		\end{equation} 
		This completes the proof.
	\end{proof}
\begin{figure}[h]
	\centering
	\includegraphics[width=0.3\textwidth,height=0.35\textwidth]{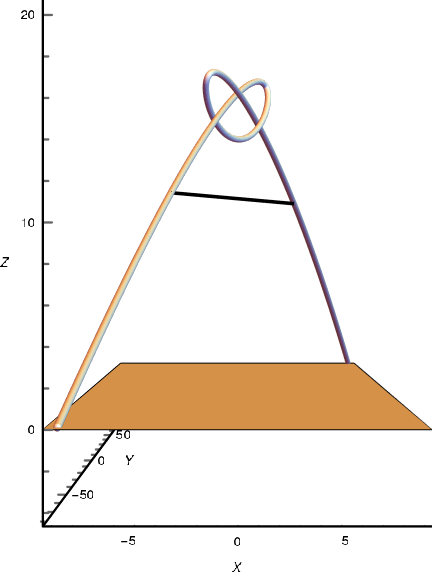}
	\caption{knotted trefoil arc with the axis of rotation}
	\label{fig: trefoil with axis}
\end{figure}	
	\begin{example}[d-twist spun trefoil]
		Start with a polynomial representation of the long trefoil knot (Figure \ref{fig: trefoil with axis}), given by
	
		\begin{align*}
			x(t)&=t^3-3t \\
			y(t)&=t^5-10t \\
			z(t)&=-t^4+4t^2+16 .
		\end{align*}
		
		Solving for $t_{1},t_{2}$ such that $x(t_{1})
		=x(t_{2})$ and $y(t_{1})
		=y(t_{2})$ we can get the crossing data and solving $z(t)=0$ will give the interval for the arc. Therefore, we have \begin{align*}
			[a,b]&=[-2.54404,2.54404],\\
			[a',b']&=[-1.946,1.946] ,\\
			P&=(f[-2.19],g[-2.19],h[-2.19]),\\
			Q&=(f[2.19],g[2.19],h[2.19]),\\
			c&=h[2.19].
		\end{align*}
		Hence, $\mathbf{k}$ is along the line joining $(-f[2.19],-g[2.19],0)$ and $(f[2.19],g[2.19],0)$.
		So,\[\mathbf{k}=\Big(\frac{f[2.19]}{\sqrt{f[2.19]^{2}+g[2.19]^{2}}},\frac{g[2.19]}{\sqrt{f[2.19]^{2}+g[2.19]^{2}}},0\Big)\]
		and \[N=\sqrt{f[2.19]^{2}+g[2.19]^{2}}.\]
		
		\begin{figure}[h]
			\centering
			\includegraphics[width=0.65\textwidth]{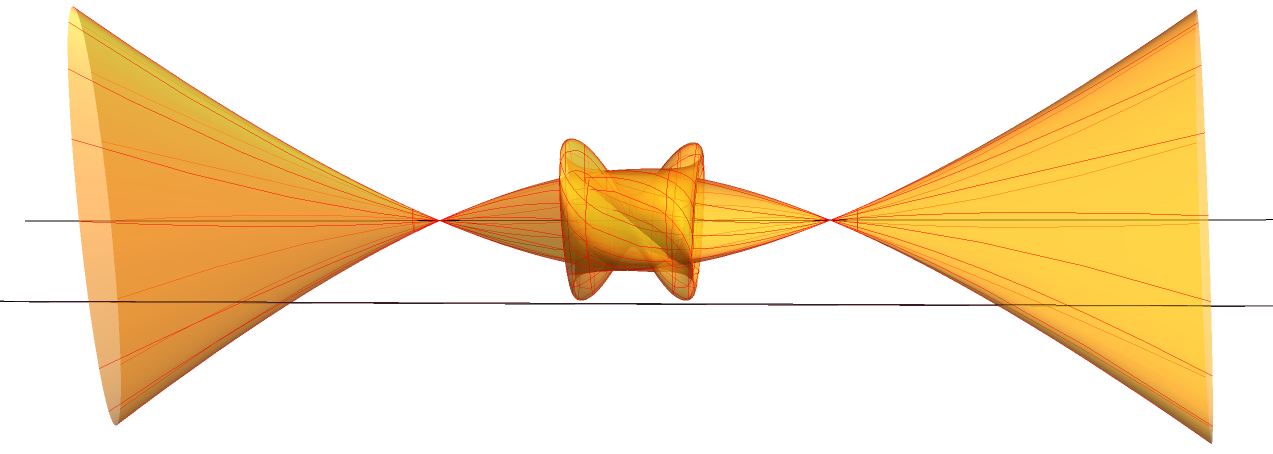}
			\caption{Rotating the long trefoil arc about $PQ$}
			\label{fig: Rotating about PQ}
		\end{figure}
		This rotation about $PQ$ is given by the parameterization (See Figure \ref{fig: Rotating about PQ}):
		\begin{equation*}
			(t,\phi) \longrightarrow  \LC \LP f'(t,\phi),g'(t,\phi),h'(t,\phi) \RP\;\middle|\;
			\begin{aligned}
				a \leq &t \leq b \\
				0 \leq &\phi < 2\pi
			\end{aligned}\, \RC 
		\end{equation*}
		where,
		{\normalsize \begin{align}
				f'(t,\phi)&=\frac{\splitdfrac{\LP f[2.19]^{2}\;+\;g[2.19]^{2}\;\cos\phi \RP\,\mathbf{f(t)}+f[2.19]\;g[2.19]\;(1-\cos\phi)\,\mathbf{g(t)}}{+(\mathbf{h(t)}-c)N\;g[2.19]\;\sin\phi}}{N^2}\\
				g'(t,\phi)&=\frac{\splitdfrac{f[2.19]\;g[2.19]\;(1-\cos\phi)\,\mathbf{f(t)}+\LP f[2.19]^{2}\;\cos\phi\;+\;g[2.19]^{2} \RP\,\mathbf{g(t)}}{+(c-\mathbf{h(t)})N\;g[2.19]\;\sin\phi}}{N^{2}}\\
				h'(t,\phi)&=\frac{-g[2.19]\;\sin\phi\;\mathbf{f(t)}+f[2.19]\sin\phi\;\mathbf{g(t)}+N
					\,\mathbf{h(t)}\cos\phi+N\,c\,(1-\cos\phi)}{N}
		\end{align}}
	
	Now, the bump function is
		\[B(t)=\frac{F(4.8-t^{2})}{F(t^{2}-3.8)+F(4.8-t^{2})}.\]
		The restricted rotation about $PQ$ (See Figure \ref{Restricted Rotation} ) is given by 
		\begin{equation*}
			(t,\phi) \longrightarrow  \LC \LP \Tilde{f}(t,\phi),\Tilde{g}(t,\phi),\Tilde{h}(t,\phi)  \RP\,\middle|\,
			\begin{aligned}
				-2.54404 \leq &t \leq 2.54404 \\
				0 \leq &\phi < 2\pi 
			\end{aligned} \RC 
		\end{equation*}
		where
		\begin{align}
			\begin{aligned}
				\Tilde{f}(t,\phi)=& \,B(t)\,f'(t,\phi)+(1-B(t))\,f(t),\\
				\Tilde{g}(t,\phi)=& \,B(t)\,g'(t,\phi)+(1-B(t))\,g(t),\\
				\Tilde{h}(t,\phi)=& \,B(t)\,h'(t,\phi)+(1-B(t))\,h(t).
			\end{aligned}.
		\end{align}
	
\begin{figure}[h]
	\centering
	\includegraphics[width=0.65\textwidth,height=0.3\textwidth]{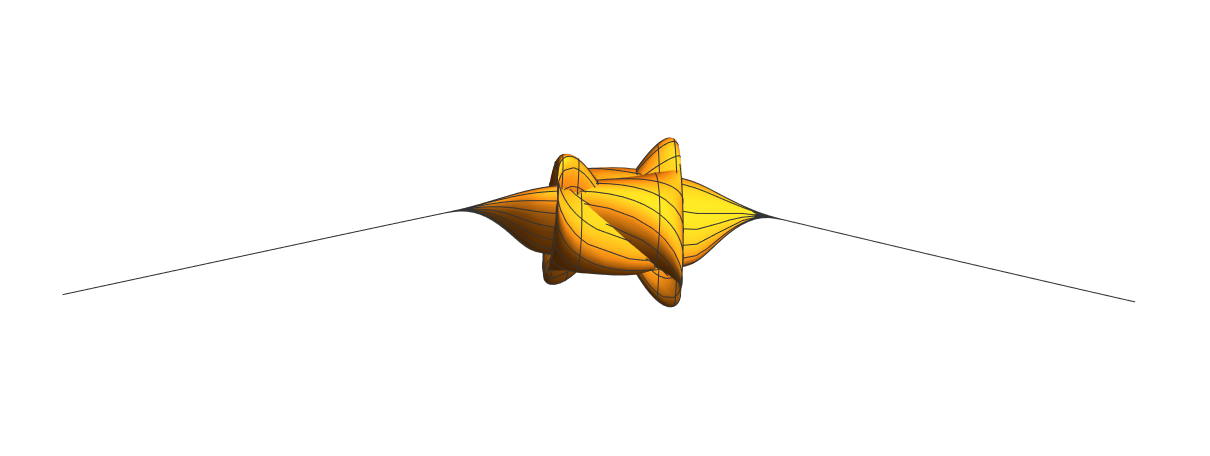}
	\caption{$( \Tilde{f}(t,\phi), \Tilde{g}(t,\phi), \Tilde{h}(t,\phi))$ in $[-2.54404,2.54404]$}
	\label{Restricted Rotation}
\end{figure}
		Hence, the parameterization for the d-twist spun trefoil knot is given by 
		{\small \begin{equation*}
				(t,\theta) \longrightarrow  \LC \LP \Tilde{f}(t,d\theta),\Tilde{g}(t,d\theta),\Tilde{h}(t,d\theta)\,C(\theta) ,\Tilde{h}(t,d\theta)\,S(\theta)  \RP\;\middle|\;
				\begin{aligned}
					-2.54404 \leq &t \leq 2.54404\\
					0 \leq &\theta < 2\pi
				\end{aligned}
				\RC .
		\end{equation*}}
		
		The {\it Mathematica} program for a projection of d-twist spun knot is provided in Section 
		\ref{mathematica}. Some of the projections of the d-twist spun trefoil on a hyperplane for the values $d=0,1,2,5,10,20$ are given below. In Figure \ref{fig: zerotwist}, it is clear that $d=0$ represents Artin's spun knot. Figure \ref{fig: zerotwistin} and Figure \ref{fig: 1twistin} illustrate the deformation of the knotted arc while rotating around $\R^{2}$ in the twist spun knot construction, which differs from Artin's spun knot construction.
	
	\begin{figure}[h!]
		\centering
		\begin{subfigure}[c]{0.4\linewidth}
			\centering
			\includegraphics[width=0.78\textwidth,height=0.45\textwidth]{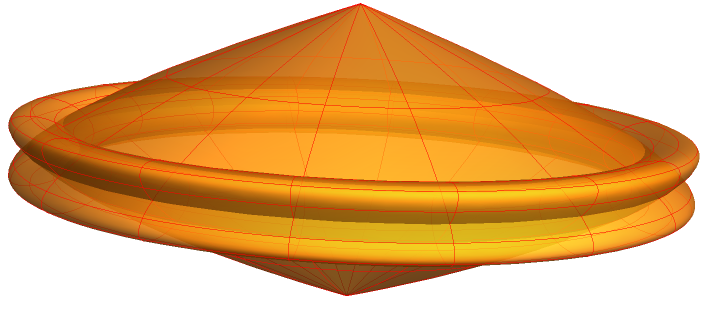}
			\caption{0-twist spun trefoil}
			\label{fig: zerotwist}
		\end{subfigure}
		\begin{subfigure}[c]{0.4\linewidth}
			\centering
			\includegraphics[width=0.78\textwidth]{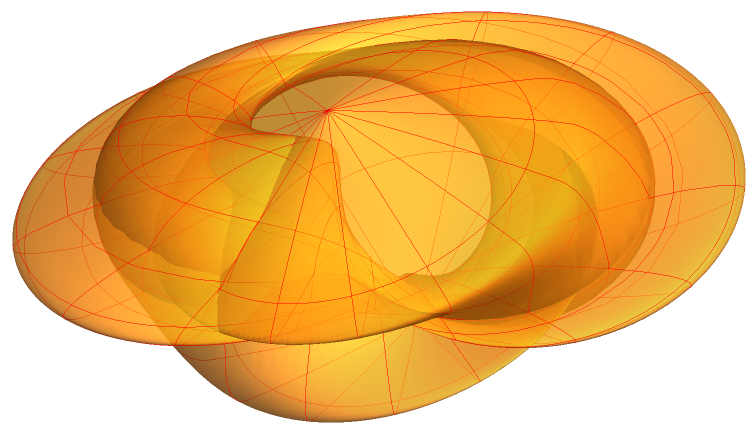}
			\caption{1-twist spun trefoil}
			\label{fig: 1twist}
		\end{subfigure}
		\newline
		\begin{subfigure}[c]{0.4\linewidth}
			\centering
			\includegraphics[width=0.75\textwidth,height=0.45\textwidth]{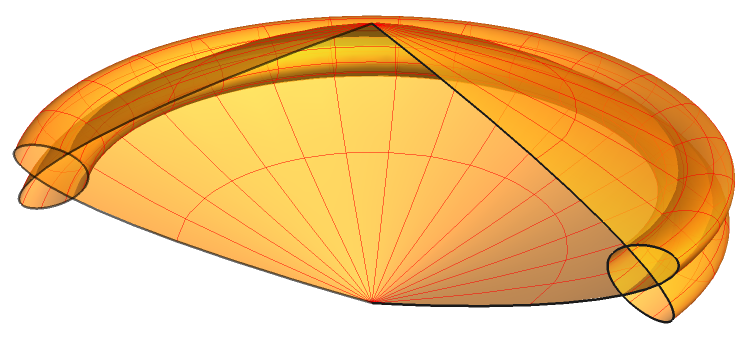}
			\caption{0-twist spun trefoil: inside view}
			\label{fig: zerotwistin}
		\end{subfigure}
		\begin{subfigure}[c]{0.45\linewidth}
			\centering
			\includegraphics[width=0.75\textwidth]{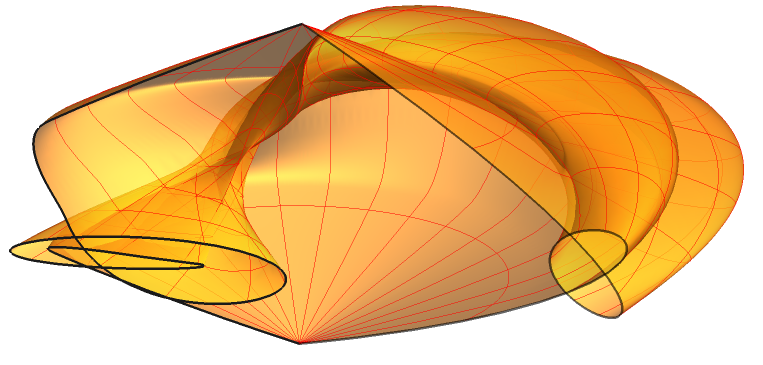}
			\caption{1-twist spun trefoil: inside view}
			\label{fig: 1twistin}
		\end{subfigure}
			\caption{d-twist spun 2-trefoil for d=0,1}
		\end{figure}
		
\begin{figure}[h!]
			\centering
		\begin{subfigure}[c]{0.45\linewidth}
			\centering
			\includegraphics[width=0.7\textwidth]{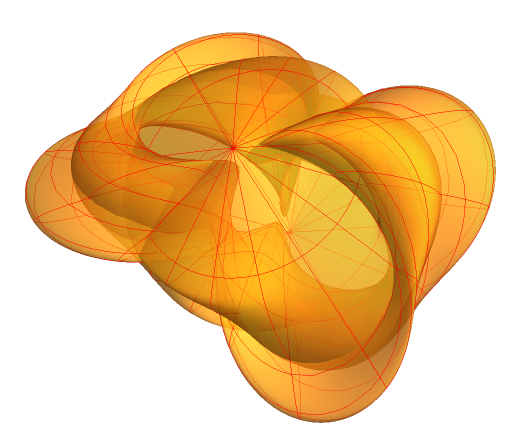}
			\caption{2-twist spun 2-trefoil}
			\label{fig: 2twist}
		\end{subfigure}
		\begin{subfigure}[c]{0.45\linewidth}
			\centering
			\includegraphics[width=0.7\textwidth]{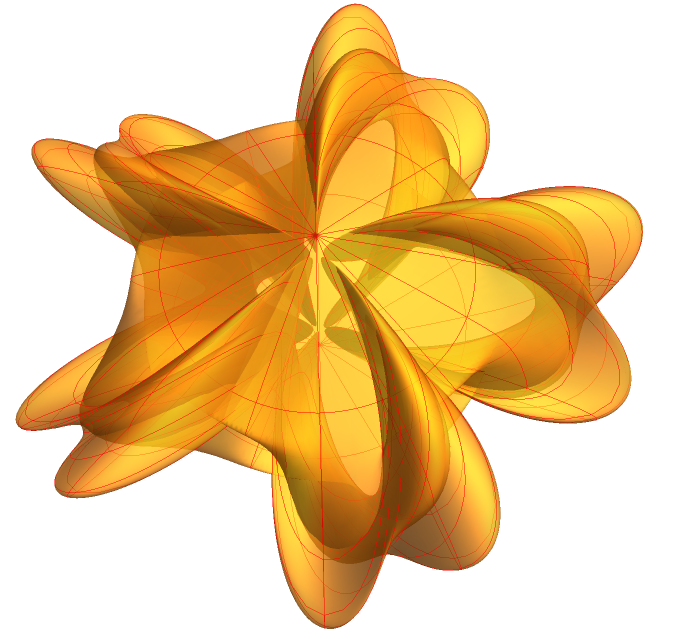}
			\caption{5-twist spun 2-trefoil}
			\label{fig: 5twist}
		\end{subfigure} 
		\newline
		\begin{subfigure}[c]{0.45\linewidth}
			\centering
			\includegraphics[width=0.7\textwidth]{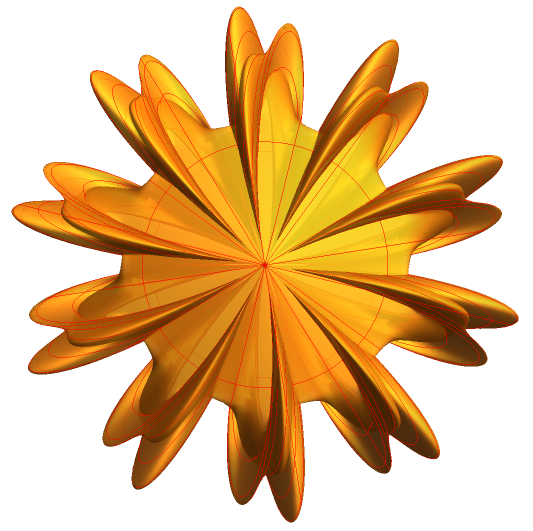}
			\caption{10-twist spun 2-trefoil}
			\label{fig: 10twist}
		\end{subfigure}
		\begin{subfigure}[c]{0.45\linewidth}
			\centering
			\includegraphics[width=0.7\textwidth]{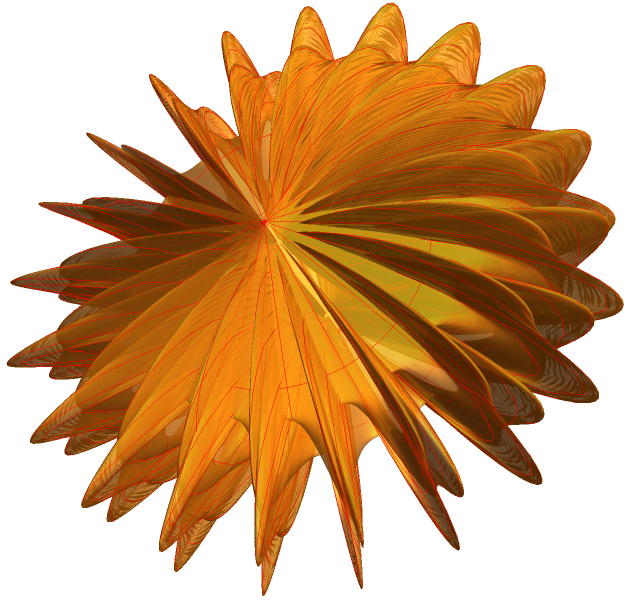}
			\caption{20-twist spun 2-trefoil}
			\label{fig: 20twist}
		\end{subfigure}
		\caption{d-twist spun 2-trefoil for d=2,5,10,20}
\end{figure}
\end{example}

\section{Conclusion:}
A more general spinning construction is described by R.A.Litherland \cite{Lit} called \say {\textit{Deform Spinning}}, where spun knots and d-twist spun knots are special cases. We would like to find a way to parameterize these $2$-knots since these are obtained by deforming classical knotted arc using some transformations in $\R^{3}$ while spinning. Finding a general method to parameterize any $2$-knot will be a goal in the next project.

\section{Acknowledgment}  
The authors are thankful to  Louis Kauffman and Visakh Narayanan for many valuable discussions. The first author would like to thank CSIR for its support.

\includepdf[width=\paperwidth,pages={1},trim={0in 0in 0in 0in},
pagecommand={\section{Mathematica Notebooks}\label{mathematica}},offset={0.5in -0.5in}]{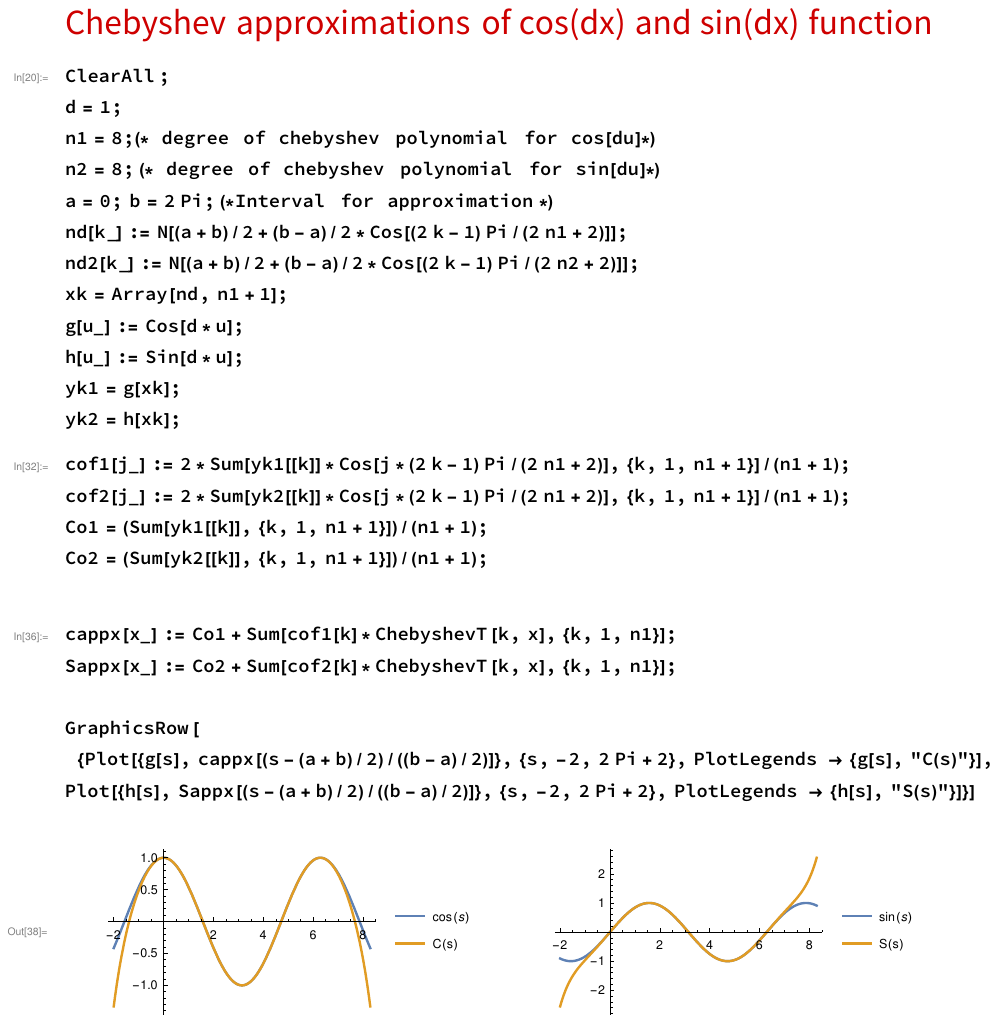}

\includepdf[pages=-,noautoscale,width=0.9\paperwidth]{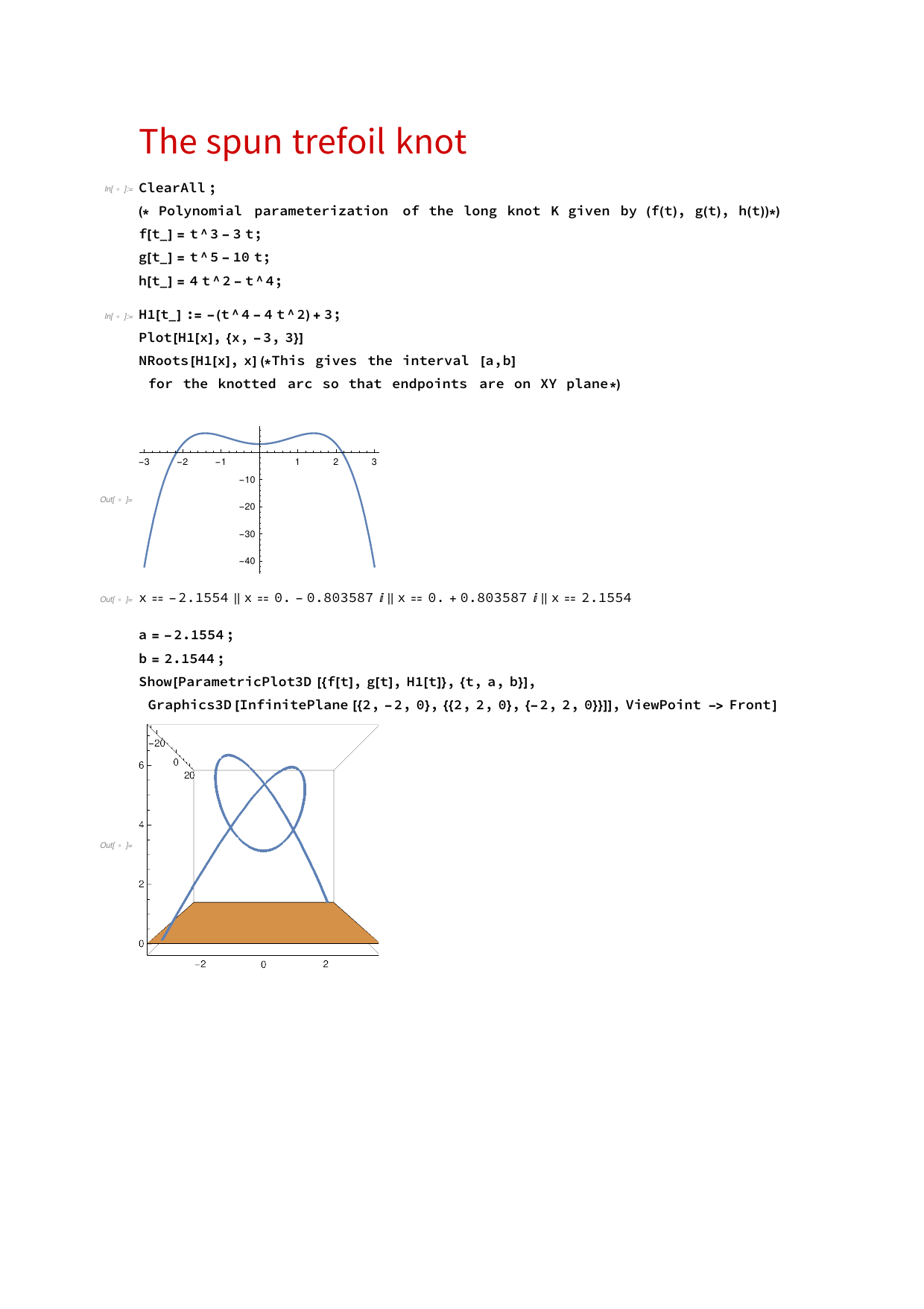}
\includepdf[pages=-,noautoscale,width=0.9\paperwidth]{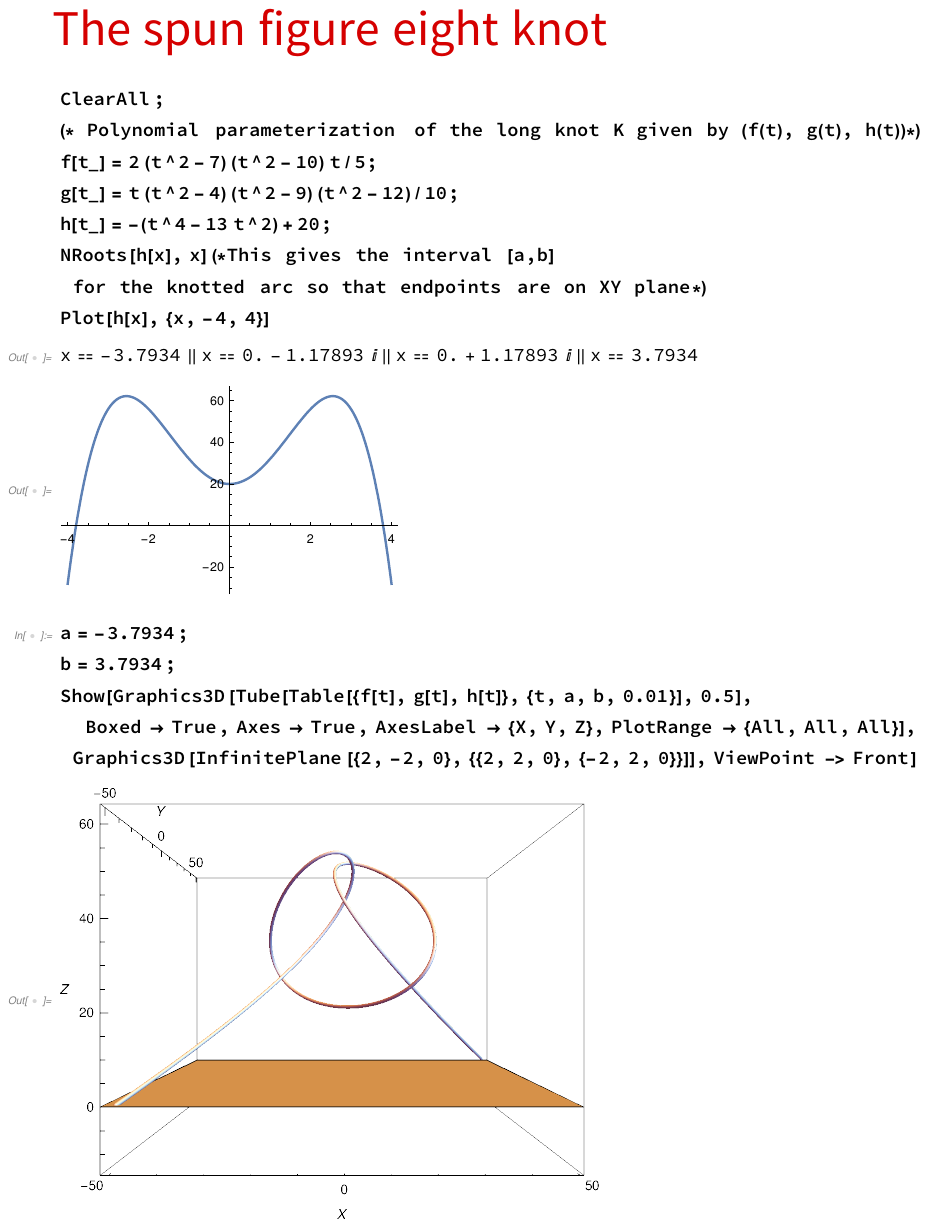}
\includepdf[pages=-,noautoscale,width=0.9\paperwidth]{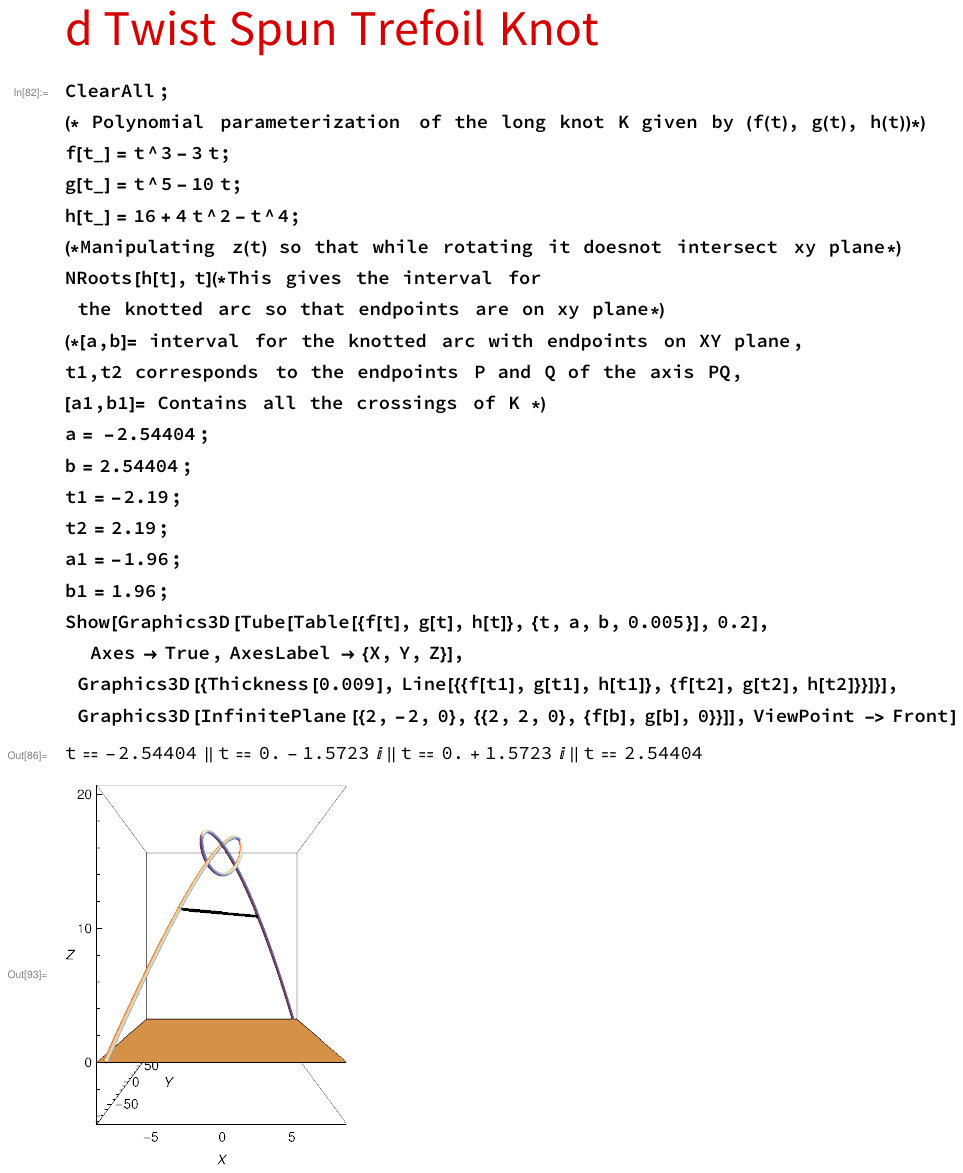}

\end{document}